\documentclass[12pt,english,refpage,intoc,bibliography=totoc,index=totoc,BCOR7.5mm,captions=tableheading]{amsart}
\usepackage{lmodern}

\usepackage[T1]{fontenc}
\usepackage[latin9]{inputenc}
\usepackage{geometry}
\usepackage{color}
\usepackage{babel}
\usepackage{textcomp}
\usepackage{enumitem}
\usepackage{amstext}
\usepackage{amsthm}
\usepackage{amssymb}
\usepackage{cancel}
\PassOptionsToPackage{normalem}{ulem}
\usepackage{ulem}
\usepackage[unicode=true,
 bookmarks=true,bookmarksnumbered=true,bookmarksopen=false,
 breaklinks=false,pdfborder={0 0 1},backref=false,colorlinks=true]
 {hyperref}
\hypersetup{pdftitle={The LyX User's Guide},
 pdfauthor={LyX Team},
 pdfsubject={LyX},
 pdfkeywords={LyX},
 linkcolor=black, citecolor=black, urlcolor=blue, filecolor=blue, pdfpagelayout=OneColumn, pdfnewwindow=true, pdfstartview=XYZ, plainpages=false}

\makeatletter
\numberwithin{equation}{section}
\numberwithin{figure}{section}
\theoremstyle{remark}
\newtheorem{note}{\protect\notename}
\theoremstyle{plain}
\newtheorem{prop}{\protect\propositionname}[section]
\theoremstyle{definition}
\newtheorem{defn}{\protect\definitionname}[section]
\theoremstyle{remark}
\newtheorem{rem}{\protect\remarkname}[section]
\theoremstyle{remark}
\newtheorem{notation}{\protect\notationname}
\theoremstyle{definition}
\newtheorem{example}{\protect\examplename}[section]
\theoremstyle{plain}
\newtheorem{lem}{\protect\lemmaname}[section]
\theoremstyle{plain}
\newtheorem{thm}{\protect\theoremname}[section]

\usepackage[figure]{hypcap}

\let\myTOC\tableofcontents
\renewcommand\tableofcontents{%
  \frontmatter
  \pdfbookmark[1]{\contentsname}{}
  \myTOC
  \mainmatter }

\usepackage{fancyhdr}
\@ifundefined{extrarowheight}
 {\usepackage{array}}{}
\makeatother

\providecommand{\definitionname}{Definition}
\providecommand{\examplename}{Example}
\providecommand{\lemmaname}{Lemma}
\providecommand{\notationname}{Notation}
\providecommand{\notename}{Note}
\providecommand{\propositionname}{Proposition}
\providecommand{\remarkname}{Remark}
\providecommand{\theoremname}{Theorem}

\begin{document}
\title{The Hydra Map and Numen Formalisms for Collatz-Type Problems}
\author{by Maxwell C. Siegel }
\curraddr{1626 Thayer Ave., Los Angeles, CA, 90024}
\email{siegelmaxwellc@ucla.edu}
\date{17 January, 2026}
\begin{abstract}
This paper details a generalization of the formalism presented in
the author's 2024 paper, \textquotedblleft The Collatz Conjecture
\& Non-Archimedean Spectral Theory: Part I \textemdash{} Arithmetic
Dynamical Systems and Non-Archimedean Value Distribution Theory\textquotedblright ,
to the case of Hydra maps on the ring of integers $\mathcal{O}_{K}$
of a global field $K$. In addition to recounting these definitions,
background material is presented for the necessary standard material
in algebraic number theory and integration and Fourier analysis with
respect to the $p$-adic Haar measure. This paper is meant to serve
as a technical manual for use of Hydra maps and numens in future research. 
\end{abstract}

\maketitle

\section{\label{sec:Document-Summary}Document Summary}

As explained in the abstract, this paper is meant to be a foundational
reference for works involving the class of generalized Collatz-type
maps I call \textbf{Hydra maps}. Its primary purpose is to fix certain
notations, definitions, conventions, and formalisms that are used
throughout my works of research on this topic.\textbf{ }Of particular
import is the aim of giving a formalism which deals with both the
archimedean and non-archimedean completions of the underlying number
field in a unified way. Future versions of this paper will hopefully
extend this unification to a fully rigorous adèlic formulation. 
\begin{note}
Readers interested only in definitions may safely skip all proofs.
Readers familiar with algebraic number theory may skip Section 3 entirely.
Sections 6 \& 7 contain the only genuinely novel analytic constructions.
The definitions for the class of maps being studied are given in Section
5.
\end{note}
\textbf{Section \ref{sec:Preliminary-Notation}} contains basic notation,
about 90\% of which is completely standard. Non-standard notations
include the author's use of $\overset{a}{\equiv}$ to indicate congruence
mod $a$ and his use of $\left[\cdot\right]_{a}$ to denote the standard
representative element set of representatives of projections mod $a$.
These notations, though minor, are crucial, as my formalism becomes
unreadable without them.

\textbf{Section \ref{sec:Global Fields}} contains a refresher on
the definitions and notations for seminorms and absolute values on
rings and fields, and their use in defining the standard concept of
a \textbf{place }as used in algebraic number theory, and using uniformizers
to extend $p$-adic values to field extensions. 

\textbf{Section \ref{sec:Dynamical-Systems}} contains standard, elementary
definitions and results from the theory of dynamical systems, with
the concept of a divergent trajectory being defined with respect to
a given absolute value on a field.

\textbf{Section \ref{sec:Hydra-Maps} }presents the definition of
\textbf{Hydra maps}, the class of generalized Collatz-type map I am
studying. This terminology is entirely non-standard, and therefore
should not be skipped. \uline{Note}: due to the rapidly evolving
nature of this subject, earlier definitions of Hydras (such as those
given in my dissertation \cite{my dissertation}) may be inconsistent
with the ones given here. The definitions given in this document will
be used as the standard going forward. For arithmetical purposes,
Hydras are defined here as maps on the ring of integers of a given
global field, however, by linear algebra, they can be defined on vector
spaces over $\mathbb{Q}$ or over a field of rational functions over
a finite field. Moreover, the definition can certainly be generalized
to Dedekind domains, though I will not do so in this version of this
paper.

\textbf{Section \ref{sec:The-Numen-of}} introduces the string formalism
(adapted from the theory of automata; see \cite{Automatic sequences})
and uses it to construct the numen, $X_{H}$, of a given Hydra map
$H$. The characterizing functional equation of $X_{H}$ is established,
and \textbf{Lemma \ref{lem:extension and measurability of X_H}} generalizes
the work in \cite{first blog paper} to show when and how $X_{H}$
may be extended from a function on $\mathbb{N}_{0}$ to a function
on the ring of $p$-adic integers. The section concludes with a statement
and sketched proof of the \textbf{Correspondence Principle}, one of
the\textbf{ }main findings of my earliest work \cite{my dissertation,first blog paper}.
This relates the integral values taken by $X_{H}$ over $\mathbb{Z}_{p}$
to the periodic points of $H$.

\textbf{Section \ref{sec:-Adic-Fourier-Analysis}} begins with a brief
review of the $p$-adic Haar probability measure and its use in $p$-adic
integration. Several basic $p$-adic integration formulae are introduced,
along with the elementary theory of Pontryagin Duality from abstract
harmonic analysis and its application to $p$-adic Fourier analysis.
Novel content in this section includes \textbf{Definition \ref{def:characteristic function}}
on page \pageref{def:characteristic function}, which gives the $\ell$-adic
characteristic function of $X_{H}$, for varying places $\ell$ of
the underlying global field $K$. The section also gives the definition
of the $\ell$-adic Wiener algebras. These are the spaces of all functions
$K_{\ell}\rightarrow\mathbb{C}$ expressible as the inverse Fourier
transform of an $L^{1}$ function $K_{\ell}\rightarrow\mathbb{C}$,
and provide a natural functional-analytic backdrop for giving a unified
treatment of Fourier analysis across the various completions of $K$.
Using the Wiener algebras, we can realize the $\ell$-adic characteristic
function of $H$ as the Fourier transform of a self-similar Borel
probability measure $d\mu_{H,\ell}$. (Equivalently, this measure
is the probability distribution associated to the $\ell$-adic valued
random variable induced by $X_{H}$.) The section concludes with some
useful formulae for working with $X_{H}$ as an $\ell$-adic valued
random variable in the case where $\ell$ is a non-trivial non-archimedean
place of $K$.

\section{\label{sec:Preliminary-Notation}Preliminary Notation}

\subsection{Notations for Sets}

For any prime $p$, we write $\mathbb{Z}_{p}$ and $\mathbb{Q}_{p}$
to denote the ring of $p$-adic integers and fields of $p$-adic rational
numbers, respectively, equipped with the standard $p$-adic absolute
value $\left|\cdot\right|_{p}$, with $\left|p\right|_{p}=1/p$. $\mathbb{Z}_{p}^{\times}$
denotes the group of multiplicatively invertible $p$-adic integers;
this is the set of all $p$-adic integers which are not congruent
to $0$ mod $p$. We write $\mathbb{C}_{p}$ to denote the field of
$p$-adic complex numbers, i.e., the metric completion of the algebraic
closure of $\mathbb{Q}_{p}$. Note that $\mathbb{C}_{p}$ is \emph{not
}spherically complete.

We write $v_{p}\left(\cdot\right)$ to denote the $p$-adic valuation,
with $\left|\cdot\right|_{p}=p^{-v_{p}\left(\cdot\right)}$ and $v_{p}\left(0\right)\overset{\textrm{def}}{=}\infty$,
where $\overset{\textrm{def}}{=}$ means ``by definition''.

For any integer $N\geq1$, we write $\mathbb{Z}/N\mathbb{Z}$ to denote
the set $\left\{ 0,\ldots,N-1\right\} $. We write $\left(\mathbb{Z}/N\mathbb{Z}\right)^{\times}$
to denote the subset of $\left\{ 0,\ldots,N-1\right\} $ consisting
of all integers co-prime to $N$. We will treat these sets as the
ring of integers mod $N$ and the multiplicative group of units mod
$N$, respectively, as the need arises.

We write $\hat{\mathbb{Z}}_{p}$ to denote $\mathbb{Z}\left[1/p\right]/\mathbb{Z}$,
the Pontryagin dual of $\mathbb{Z}_{p}$. We identify $\hat{\mathbb{Z}}_{p}$
with the group of rational numbers in $\left[0,1\right)$ whose denominators
are powers of $p$. Viewing $\hat{\mathbb{Z}}_{p}$ as a subset of
$\mathbb{Q}_{p}$, for any $t\in\hat{\mathbb{Z}}_{p}\backslash\left\{ 0\right\} $,
observe that $t$ can be written in irreducible form as $k/p^{n}$
for some integers $n\geq1$ and $k\in\left(\mathbb{Z}/p\mathbb{Z}\right)^{\times}$.
Consequently, the $p$-adic absolute value of $t$ is $\left|t\right|_{p}=p^{n}$,
and the numerator of $t$ is given by $t\left|t\right|_{p}=k$, with
both values being $0$ when $t=0$. We write $u_{p}:\mathbb{Q}_{p}\rightarrow\mathbb{Z}_{p}^{\times}\cup\left\{ 0\right\} $
to denote the function $u_{p}\left(\mathfrak{y}\right)\overset{\textrm{def}}{=}\mathfrak{y}\left|\mathfrak{y}\right|_{p}$.

For any real number $x$, we write $\mathbb{N}_{x}$ to denote the
set of all integers $\geq x$ (thus, $\mathbb{N}_{0}=\left\{ 0,1,2,\ldots\right\} $;
$\mathbb{N}_{1}=\left\{ 1,2,\ldots\right\} $, etc.). We write $\mathbb{Z}_{p}^{\prime}$
to denote $\mathbb{Z}_{p}\backslash\mathbb{N}_{0}$; note that this
is the set of all $p$-adic integers with infinitely many non-zero
$p$-adic digits. $\mathbb{P}$, meanwhile, denotes the set of prime
numbers: $2,3,5,7,\ldots$.

Also, as the reader may have noticed, for aesthetic reasons, we write
$p$-adic variables in lower-case $\mathfrak{fraktur}$ font.

\subsection{Notation for Congruences}

For any $\mathfrak{z}\in\mathbb{Z}_{p}$ and any integer $n\geq0$,
we write $\left[\mathfrak{z}\right]_{p^{n}}$ to denote the unique
integer in $\mathbb{Z}/p^{n}\mathbb{Z}$ which is congruent to $\mathfrak{z}$
modulo $p^{n}$. For terminological purposes, we follow Yvette Amice
in referring to the series representation $\sum_{n=-n_{0}}^{\infty}c_{n}p^{n}$
of a $p$-adic number $\mathfrak{y}\in\mathbb{Q}_{p}$ (be it rational
or integral) as the \textbf{Hensel series} of $\mathfrak{y}$.

We adopt the convention of viewing $p$-adic integers as power series
(and, likewise, $p$-adic rational numbers as Laurent series) in the
variable $p$, and, importantly we write $p$-adic digits from \emph{left
to right}, in order of increasing powers of $p$. Thus, in $2$-adic
notation, we would write:
\begin{align*}
1 & =\centerdot_{2}1\\
2 & =\centerdot_{2}01\\
3 & =\centerdot_{2}11\\
4 & =\centerdot_{2}001\\
 & \vdots
\end{align*}
and so on. Here, we use $\centerdot_{2}$ as the ``$2$-adic decimal
point''. Throughout this series, all $2$-adic numbers written in
base $2$ will have a $\centerdot_{2}$ in them, and, more generally,
all $p$-adic numbers written in base $p$ will have a $\centerdot_{p}$
in them. For example:
\begin{equation}
11\centerdot_{2}011=\frac{1}{2^{2}}+\frac{1}{2}+2^{1}+2^{2}=6+\frac{3}{4}
\end{equation}

\begin{equation}
121\centerdot_{3}01=\frac{1}{3^{3}}+\frac{2}{3^{2}}+\frac{1}{3}+3^{1}=\frac{16}{27}+3
\end{equation}

We employ a non-standard notation for congruences:
\begin{align}
x & \overset{a}{\equiv}y\nonumber \\
 & \Updownarrow\\
x & =y\mod a\nonumber 
\end{align}
To give some examples, for $\mathfrak{z}\in\mathbb{Z}_{p}$, $\mathfrak{z}\overset{p^{n}}{\equiv}k$
means ``$\mathfrak{z}$ is congruent to $k$ mod $p^{n}$''; i.e.,
$\mathfrak{z}\in k+p^{n}\mathbb{Z}_{p}$. Given $\mathfrak{x},\mathfrak{y}\in\mathbb{Q}_{p}$,
we write $\mathfrak{x}\overset{1}{\equiv}\mathfrak{y}$ to mean $\mathfrak{x}-\mathfrak{y}=0\mod\mathbb{Z}_{p}$.
In particular, note that:
\begin{itemize}
\item $\mathfrak{y}\overset{1}{\equiv}\left\{ \mathfrak{y}\right\} _{p}$
for all $\mathfrak{y}\in\mathbb{Q}_{p}$;
\item $s,t\in\hat{\mathbb{Z}}_{p}$ denote the same element of $\hat{\mathbb{Z}}_{p}$
if and only if $s\overset{1}{\equiv}t$;
\item $\mathfrak{z}\overset{1}{\equiv}k$ is true for all $\mathfrak{z}\in\mathbb{Z}_{p}$
and all $k\in\mathbb{Z}$.
\end{itemize}
For a rational integer $n$, $n\overset{2}{\equiv}0$ means $n$ is
even, while $n\overset{2}{\equiv}1$ means $n$ is odd. Naturally,
we write $x\overset{a}{\cancel{\equiv}}y$ to mean ``$x$ is \emph{not
}congruent to $y$ mod $a$''.

\subsection{Notation for Convergence}

Let $K$ be a topological field. Then, in equations with limits or
infinite sums, we write $\overset{K}{=}$ to mean that the convergence
occurs with respect to the topology of $K$. For example, $\overset{\mathbb{Q}_{p}}{=}$
to indicate that the convergence is with respect to the topology of
$\mathbb{Q}_{p}$. Similarly, $\overset{\mathbb{R}}{=}$ means convergence
with respect to the topology of $\mathbb{R}$. If $K=\overline{\mathbb{Q}}$
or any subfield thereof ($\mathbb{Q}$, $\mathbb{Q}\left(\sqrt{2}\right)$,
etc.), we always equip $K$ with the discrete topology. So, $\overset{\mathbb{Q}}{=}$
and $\overset{\overline{\mathbb{Q}}}{=}$ denote convergence with
respect to the discrete topologies on $\mathbb{Q}$ and $\overline{\mathbb{Q}}$,
respectively. Note then that a sequence of rational numbers $\left\{ a_{n}\right\} _{n\geq0}$
satisfies $\lim_{n\rightarrow\infty}a_{n}\overset{\mathbb{Q}}{=}0$
if and only if $a_{n}=0$ for all sufficiently large $n$.

Also, given a topological ring $R$, we write $\overset{R}{=}$ to
mean that convergence occurs in the topology of $K$, where $K$ is
the field of fractions of $R$, and that the terms of the sequence
being limited or the partial sums of the series being limited are
also elements of $R$. Thus, for a function $f:\mathbb{Z}_{2}\rightarrow\mathbb{Z}_{q}$,
we write:
\begin{equation}
\lim_{n\rightarrow\infty}f\left(\left[\mathfrak{z}\right]_{2^{n}}\right)\overset{\mathbb{Z}_{q}}{=}f\left(\mathfrak{z}\right),\textrm{ }\forall\mathfrak{z}\in\mathbb{Z}_{2}
\end{equation}
to mean that, for each $\mathfrak{z}\in\mathbb{Z}_{2}$, the sequence
$\left\{ f\left(\left[\mathfrak{z}\right]_{2^{n}}\right)\right\} _{n\geq0}$
consists of elements of $\mathbb{Z}_{q}$ that converge in $\mathbb{Z}_{q}$
to $f\left(\mathfrak{z}\right)$ as $n\rightarrow\infty$.

\subsection{Preliminaries from $p$-Adic Analysis}

For a refresher on the topic, the reader can refer to Robert's book
\cite{Robert's Book}, Schikhof's book \cite{Ultrametric Calculus},
or the author's doctoral dissertation \cite{my dissertation}. For
our immediate purposes, we need only a few results from the theory
of ultrametric analysis:
\begin{prop}
The $p$-adic absolute value $\left|\cdot\right|_{p}$ satisfies the
\textbf{ultrametric inequality}:
\begin{equation}
\left|\mathfrak{a}+\mathfrak{b}\right|_{p}\leq\max\left\{ \left|\mathfrak{a}\right|_{p},\left|\mathfrak{b}\right|_{p}\right\} ,\textrm{ }\forall\mathfrak{a},\mathfrak{b}\in\mathbb{Q}_{p}
\end{equation}
with equality if and only if $\left|\mathfrak{a}\right|_{p}\neq\left|\mathfrak{b}\right|_{p}$.
\end{prop}
\begin{prop}
\label{prop:Non-archimedean series converge iff their terms decay to zero}Let
$\left(\mathbb{K},\left|\cdot\right|_{\mathbb{K}}\right)$ be a metrically
complete non-archimedean field, and let $\left\{ a_{n}\right\} _{n\geq0}$
be a sequence in $\mathbb{K}$. Then, the series $\sum_{n=0}^{\infty}a_{n}$
converges in $\mathbb{K}$ if and only if $\lim_{n\rightarrow\infty}\left|a_{n}\right|_{\mathbb{K}}\overset{\mathbb{R}}{=}0$.
\end{prop}
\begin{prop}
Let $p$ be a prime number. Then, a $p$-adic number $\mathfrak{y}\in\mathbb{Q}_{p}$
is a rational number (i.e., an element of $\mathbb{Q}$) if and only
if the sequence of $\mathfrak{y}$'s $p$-adic digits is eventually
periodic.
\end{prop}

\section{\label{sec:Global Fields}Global Fields, Absolute Values, \& Completions}

By a \textbf{base field} $\mathbb{F}$, we mean either the field of
rational numbers ($\mathbb{Q}$), or the field $\mathbb{F}_{q_{\mathbb{F}}}\left(T\right)$
of formal rational functions of $T$ with coefficients in $\mathbb{F}_{q_{\mathbb{F}}}$,
the finite field of characteristic ($\textrm{char}$) $p_{\mathbb{F}}$
and order $q_{\mathbb{F}}=p_{\mathbb{F}}^{n_{\mathbb{F}}}$, where
$p_{\mathbb{F}}$ is a prime number and $n_{\mathbb{F}}$ is a positive
integer. The\textbf{ ring of integers of $\mathbb{F}$}, denoted $\mathcal{O}_{\mathbb{F}}$,
is either $\mathbb{Z}$ (if $\mathbb{F}=\mathbb{Q}$) or $\mathbb{F}_{q_{\mathbb{F}}}\left[T\right]$,
the ring of formal polynomials in $T$ with coefficients in $\mathbb{F}_{q_{\mathbb{F}}}$
(if $\mathbb{F}=\mathbb{F}_{q_{\mathbb{F}}}\left(T\right)$). By a
\textbf{global field }$K$, we mean a finite-degree extension of $\mathbb{F}$.
The \textbf{ring of integers }of $K$ is denoted $\mathcal{O}_{K}$.
We write $d$ (or $d_{K}$, for emphasis) to denote the degree of
$K$ over $\mathbb{F}$.

Given a prime ideal $\mathfrak{Q}$ of \emph{$\mathcal{O}_{K}$},
recall that there exists a unique prime ideal $\mathfrak{P}$ of $\mathcal{O}_{\mathbb{F}}$
so that 
\begin{equation}
\mathfrak{Q}\cap\mathcal{O}_{\mathbb{F}}=\mathfrak{P}
\end{equation}
and we then say that $\mathfrak{Q}$ \textbf{lies over $\mathfrak{P}$},
and that $\mathfrak{P}$ \textbf{lies under $\mathfrak{Q}$}. 
\begin{defn}
A \textbf{seminorm }on $R$ is a function $\left\Vert \cdot\right\Vert :R\rightarrow\left[0,\infty\right)$
so that:
\begin{align}
\left\Vert 0_{R}\right\Vert  & =0\\
\left\Vert 1_{R}\right\Vert  & =1\\
\left\Vert r+s\right\Vert  & \leq\left\Vert r\right\Vert +\left\Vert s\right\Vert ,\textrm{ }\forall r,s\in R\\
\left\Vert rs\right\Vert  & \leq\left\Vert r\right\Vert \left\Vert s\right\Vert ,\textrm{ }\forall r,s\in R
\end{align}
If, in addition, $\left\Vert \cdot\right\Vert $ satisfies the ultrametric
inequality:
\begin{equation}
\left\Vert f+g\right\Vert \leq\max\left\{ \left\Vert f\right\Vert ,\left\Vert g\right\Vert \right\} 
\end{equation}
with equality whenever $\left\Vert f\right\Vert \neq\left\Vert g\right\Vert $,
we say $\left\Vert \cdot\right\Vert $ is \textbf{non-archimedean}.
We say $\left\Vert \cdot\right\Vert $ is \textbf{multiplicative }if
$\left\Vert rs\right\Vert =\left\Vert r\right\Vert \left\Vert s\right\Vert $
for all $r,s\in R$.

We call $\left\Vert \cdot\right\Vert $ a \textbf{ring} \textbf{norm}
if $\left\Vert r\right\Vert =0$ if and only if $r=0_{R}$. A \textbf{ring
absolute value }is a ring norm which is multiplicative.

A \textbf{valued ring} is a ring with an absolute value. A \textbf{normed
ring }is a ring with a norm. A \textbf{completed valued ring }is a
valued ring so that the metric induced by the absolute value is complete.
A \textbf{Banach ring }is a normed ring so that the metric induced
by the norm is complete. A \textbf{valued field }is a valued ring
which is also a field.

A \textbf{normed ($R$-)algebra }is an $R$-algebra $\mathcal{A}$
which, as a ring, is a normed ring, and which satisfies $\left\Vert r\alpha\right\Vert =\left\Vert r\right\Vert \left\Vert \alpha\right\Vert $
for all $r\in R$ and all $\alpha\in\mathcal{A}$. A \textbf{Banach
($R$-)algebra }is a normed algebra so that the metric induced by
the norm is complete.

By the \textbf{quality }of an absolute value, field, norm, normed,
space, Banach algebra, or the like, we refer to the property of whether
the thing under consideration is archimedean or non-archimedean. Thus,
we can have a seminorm of archimedean quality, a valued ring of non-archimedean
quality, and so on.

Finally, two absolute values $\left|\cdot\right|_{a}$ and $\left|\cdot\right|_{b}$
are said to be \textbf{equivalent }if there is a real number $c>0$
so that:
\begin{equation}
\left|r\right|_{b}=\left|r\right|_{a}^{c},\textrm{ }\forall r\in R
\end{equation}
An absolute value is \textbf{trivial }if $\left|r\right|=1$ for all
$r\in R\backslash\left\{ 0\right\} $. 
\end{defn}
\begin{defn}
A \textbf{place }of a ring or field is an equivalence class of equivalent
absolute values. Given a place $\ell$ of the field $K$, we then
write $K_{\ell}$ to denote the metric completion of $K$ with respect
to any absolute value in $\ell$; up to isomorphism, $K_{\ell}$ is
independent of the choice of the representative absolute value in
$\ell$. Places whose absolute values are non-archimedean are called
\textbf{finite places}; places whose absolute values are archimedean
are called \textbf{infinite places}.We write $\textrm{Pl}\left(K\right)$
to denote the set of all places of $K$, and write $\textrm{Pl}^{\times}\left(K\right)$
to denote the set of all non-trivial places of $K$. We write $\textrm{Pl}^{\infty}\left(K\right)$
to denote the set of all non-trivial archimedean places of $K$, and
write $\textrm{Pl}^{0}\left(K\right)$ to denote the set of all non-trivial
non-archimedean places of $K$. Recall that non-archimedean places
are in a bijective correspondence with the prime ideals of $\mathcal{O}_{K}$,
while archimedean places are in a bijective correspondence with the
embeddings of $K$ in $\mathbb{C}$. If $\textrm{char}\left(K\right)>0$,
$K$ has no archimedean places.

Places will be denoted by $\ell$. $K_{\ell}$ denote the completion
of $K$ with respect to absolute values in $\ell$, respectively.
The symbol $\infty$ will be used as a subscript to indicate an archimedean
place. Thus, just as $\mathbb{Q}_{p}$ is the $p$-adic completion
of $\mathbb{Q}$, $\mathbb{Q}_{\infty}$ is the archimedean completion
of $\mathbb{Q}$, otherwise known as $\mathbb{R}$. In an abuse of
notation, given a prime ideal $\mathfrak{P}\subset\mathcal{O}_{K}$,
we will sometimes write $\left|\cdot\right|_{\mathfrak{P}}$ to denote
an absolute value of the corresponding place. Given a non-archimedean
place $\ell$ of $K$, we write $v_{\ell}$ to denote the associated
valuation (note: \emph{valuation}, not absolute value). If $\mathfrak{P}$
is the prime ideal associated to $\ell$, we also write $v_{\mathfrak{P}}$
to denote $v_{\ell}$.

Given places $\ell$ of $K$ and $q$ of $\mathbb{F}$, recall that
we say $\ell$ \textbf{lies over }$q$ whenever $\ell$ and $q$ are
of the same quality, and whenever:

\textbullet{} If $\ell$ and $q$ are non-archimedean, the prime ideal
of $\mathcal{O}_{K}$ associated to $\ell$ lies over the prime ideal
of $\mathcal{O}_{\mathbb{F}}$ associated to $q$;

\textbullet{} If $\ell$ and $q$ are archimedean, $q$ is the standard
archimedean absolute value on $\mathbb{Q}$.
\end{defn}
Recall that the completion of a global field with respect to a non-trivial
place is called a \textbf{local field}. If the place $\ell$ is non-archimedean,
$\mathcal{O}_{K_{\ell}}$, the ring of integers of $K_{\ell}$, is
the completion of $\mathcal{O}_{K}$ with respect to $\ell$, and
is a \textbf{local ring}. In particular, there exists an element $\pi_{\ell}\in\mathcal{O}_{K_{\ell}}$
called the \textbf{($\ell$-adic)} \textbf{uniformizer} so that $\pi_{\ell}\mathcal{O}_{K_{\ell}}$
is the unique maximal ideal of $\mathcal{O}_{K_{\ell}}$. If $\mathfrak{P}$
is a prime ideal of $\mathcal{O}_{\mathbb{F}}$ and $p$ is a prime
number so that a uniformizer $\pi_{\mathfrak{P}}$ of $\mathcal{O}_{\mathbb{F}_{\mathfrak{P}}}$
(the ring of integers of the $\mathfrak{P}$-adic completion of $\mathbb{F}$)
satisfies $\left|\pi_{\mathfrak{P}}\right|_{\mathfrak{P}}=p^{-1}$,
recall there is then there a unique integer $e_{K,\ell}\geq1$, called
the \textbf{ramification index of $\ell$ over $\mathfrak{P}$, }so
that:
\begin{equation}
\left|\pi_{\ell}\right|_{\ell}=p^{-1/e_{K,\ell}}
\end{equation}
Moreover, we have that $\left|\cdot\right|_{\mathfrak{P}}$ extends
uniquely to $K_{\ell}$, with:
\begin{equation}
\left|\pi_{\ell}\right|_{\mathfrak{P}}=p^{-1/e_{K,\ell}}
\end{equation}

\begin{defn}
We write $\textrm{End}_{\mathbb{F}}\left(K\right)$ to denote the
\textbf{$\mathbb{F}$-algebra of all $\mathbb{F}$-linear maps $K\rightarrow K$};
we can view these as $d\times d$ matrices with entries in $\mathbb{F}$.
We write $\textrm{End}_{\mathbb{F}}^{\times}\left(K\right)$ to denote
the \textbf{group of invertible $\mathbb{F}$-linear maps $K\rightarrow K$};
this is a group under map composition.
\end{defn}
\begin{defn}
Given a place $\ell$ of $K$ and an $r\in\textrm{End}_{\mathbb{F}}\left(K\right)$,
we write $\left\Vert r\right\Vert _{\ell}$ to denote the \textbf{$\ell$-adic
operator norm }of $r$, defined by:
\begin{equation}
\left\Vert r\right\Vert _{\ell}\overset{\textrm{def}}{=}\sup_{\begin{array}{c}
z\in K\\
\left|z\right|_{\ell}\leq1
\end{array}}\frac{\left|rz\right|_{\ell}}{\left|z\right|_{\ell}}
\end{equation}
so that $\left|rz\right|_{\ell}\leq\left\Vert r\right\Vert _{\ell}\left|z\right|_{\ell}$
for all $z\in K$. 
\end{defn}
\begin{rem}
Since $K$ is dense in the metric completion $K_{\ell}$, observe
that any $r\in\textrm{End}_{\mathbb{F}}K$ uniquely extends to a map
$K_{\ell}\rightarrow K_{\ell}$ possessing the same $\ell$-adic norm
as $r$. 
\end{rem}
\begin{defn}
Given a non-zero prime ideal $\mathfrak{P}\subseteq\mathcal{O}_{K}$,
we say $u\in\textrm{End}_{\mathbb{F}}^{\times}\left(K\right)$ is
a \textbf{$\mathfrak{P}$-adic isometry} whenever:
\begin{equation}
\left|uz\right|_{\mathfrak{P}}=\left|z\right|_{\mathfrak{P}},\textrm{ }\forall z\in K
\end{equation}
\end{defn}
\begin{rem}
Since $K$ is dense in the metric completion $K_{\mathfrak{P}}$,
observe that any $u\in\textrm{End}_{\mathbb{F}}^{\times}\left(K\right)$
which is a $\mathfrak{P}$-adic isometry then uniquely extends to
a map $u:K_{\mathfrak{P}}\rightarrow K_{\mathfrak{P}}$, which is
an isometry of $K_{\mathfrak{P}}$. 
\end{rem}
\begin{prop}[$\mu u$ Factorization (Non-Archimedean)]
\label{prop:NA mu u factorization}Let $r\in\textrm{End}_{\mathbb{F}}^{\times}\left(K\right)$,
and let $P=\left\{ \mathfrak{P}_{1},\ldots,\mathfrak{P}_{n}\right\} $
be a collection of finitely many non-zero prime ideals of $\mathcal{O}_{K}$,
with associated valuations $v_{\mathfrak{P}_{1}},\ldots,v_{\mathfrak{P}_{n}}$
and normalized non-archimedean absolute values $\left|\cdot\right|_{\mathfrak{P}_{1}},\ldots,\left|\cdot\right|_{\mathfrak{P}_{n}}$.
Then, there exist a scalar $\mu\in K^{\times}$ and a $u\in\textrm{End}_{\mathbb{F}}^{\times}\left(K\right)$
so that $r$ decomposes in $\textrm{End}_{\mathbb{F}}^{\times}\left(K\right)$
as:
\begin{equation}
r=\varpi u\label{eq:mu u factorization}
\end{equation}
where we identify $\varpi$ with the map $z\mapsto\varpi z$ in $\textrm{End}_{\mathbb{F}}^{\times}\left(K\right)$,
and where $u$ is a $\mathfrak{P}_{m}$-adic isometry for all $m\in\left\{ 1,\ldots,n\right\} $
\begin{equation}
\left|uz\right|_{\mathfrak{P}_{m}}=\left|z\right|_{\mathfrak{P}_{m}},\textrm{ }\forall z\in K\label{eq:lambda adic isometry condition}
\end{equation}
Moreover, as constructed, we have that:
\begin{align}
\left|rz\right|_{\mathfrak{P}_{m}} & =\left|\varpi z\right|_{\mathfrak{P}_{m}}
\end{align}
for all $z\in K$ and all $m\in\left\{ 1,\ldots,n\right\} $.

We call (\ref{eq:mu u factorization}) the \textbf{decomposition of
$r$ over $P$}. Given an ideal $\Lambda$ of $\mathcal{O}_{K}$,
we call (\ref{eq:mu u factorization}) the \textbf{decomposition of
$r$ over $\Lambda$} when $P$ is the set of prime ideals in the
factorization of $\Lambda$.
\end{prop}

\section{\label{sec:Dynamical-Systems}Dynamical Systems}
\begin{note}
\cite{Dynamical Systems} is an easy introduction to the subject.
\end{note}
Fix a set $X\rightarrow X$ and a map $T:X\rightarrow X$. 
\begin{notation}
For any integer $k\geq0$, we write $T^{\circ k}$ to denote the composition
of $k$ copies of $T$, with $T^{\circ0}$ being defined as the identity
map.
\end{notation}
\begin{defn}
Recall that a \textbf{periodic point }of $T$ is an $x\in X$ so that
$T^{\circ k}\left(x\right)=x$ for some $k\geq1$. We say $x$ is
\textbf{pre-periodic }if there exists a $k\geq1$ so that $T^{\circ k}\left(x\right)$
is a periodic point of $T$. We say $x$ is a \textbf{strictly pre-periodic
}point if $x$ is a pre-periodic point of $T$ which is \emph{not
}a periodic point of $T$. The \textbf{forward orbit }(sometimes called
the \textbf{trajectory})\textbf{ }of $x$ under $T$ is the sequence:
\begin{equation}
\left\{ T^{\circ k}\left(x\right)\right\} _{k\geq0}
\end{equation}
A \textbf{cycle }$\Omega$ of $T$ is the forward orbit of a periodic
point of $T$. 

If $X$ is a valued field with absolute value $\left|\cdot\right|_{X}$,
we say $x\in X$ is a \textbf{divergent point }of $T$ if : 
\begin{equation}
\lim_{k\rightarrow\infty}\left|T^{\circ k}\left(x\right)\right|_{X}\overset{\mathbb{R}}{=}+\infty
\end{equation}
We use the phrase ``a \textbf{divergent trajectory }of $T$'' to
refer to the forward orbit of a divergent point (which may or may
not exist, depending on $T$).
\end{defn}
Next, we have orbit classes and their associated definitions
\begin{defn}
Recall that an \textbf{orbit class }of $T$ is a non-empty set $V\subseteq X$
so that the pre-image of $V$ under $T$ is equal to $V$: $T^{-1}\left(V\right)=V$.

An orbit class $V$ of $T$ is said to \textbf{irreducible }if it
cannot be written as $V=V_{1}\cup V_{2}$, where $V_{1}$ and $V_{2}$
are orbit classes of $T$ with $V_{1}\cap V_{2}=\varnothing$.

A \textbf{backward orbit }of some $x\in X$ is a sequence $\left\{ x_{n}\right\} _{n\geq1}$
in $X$ so that $x_{n+1}\in T^{-1}\left(\left\{ x_{n}\right\} \right)$
for all $n\geq0$. Note that if $T$ is not injective, a given $x_{0}\in X$
can have more than one backward orbit.

We say an orbit class $V$ is \textbf{periodic }if it is irreducible
and contains a cycle. We say an orbit class $V$ is \textbf{divergent
}if it is irreducible and contains a divergent point.
\end{defn}
Lastly, we recall two elementary results from the theory of discrete
dynamical systems. For proofs, see \cite{first blog paper}.
\begin{prop}
Let $T:X\rightarrow X$ be any map. Then, the collection of irreducible
orbit classes of $T$ form a partition of $X$ into pair-wise disjoint
sets. These sets are the equivalence classes of $X$ under the relation
$\sim$ defined by:
\begin{equation}
x\sim y\Leftrightarrow\exists m,n\geq0:T^{\circ m}\left(x\right)=T^{\circ n}\left(y\right)
\end{equation}
In particular, if $X$ is countable, then $X$ consists of at most
countably many irreducible orbit classes of $T$.
\end{prop}
\begin{prop}
\label{prop:Classification of T_3's dynamics}Let $X$ be a valued
field with absolute value $\left|\cdot\right|_{X}$, and let $Y$
be a discrete subset of $X$, in the sense that there is a $\delta>0$
so that:
\begin{equation}
\inf_{\begin{array}{c}
x,y\in Y\\
x\neq y
\end{array}}\left|x-y\right|_{X}\geq\delta
\end{equation}
Then, for any map $T:Y\rightarrow Y$ be any map, every irreducible
orbit class of $T$ in $Y$ is either periodic or divergent. Consequently,
any $x\in Y$ is either a divergent point of $T$ or\emph{ }a pre-periodic\textbf{
}point of $T$.
\end{prop}
Proof: Let $V$ be an irreducible orbit class of $T$ in $Y$.

If $V$ contains a cycle $\Omega$ of $T$, then, given any $x\in V$,
$x\sim\omega$ for some periodic point $\omega\in\Omega$. So, $x$
is eventually iterated to some element of $\omega$'s trajectory.
Since \emph{every} element of $\omega$'s trajectory is in $\Omega$,
$x$ is eventually iterated into $\Omega$, which makes $x$ into
a pre-periodic point of $T$.

So, suppose that $V$ does \emph{not} contain a cycle, and, again,
let $x\in V$ be arbitrary. Since $T^{-1}\left(V\right)=V$, applying
$T$ on both sides yields:
\begin{equation}
T\left(V\right)=T\left(T^{-1}\left(V\right)\right)\subseteq V
\end{equation}
So, everything in $x$'s trajectory will also be in $V$. Thus, if
$x$ was pre-periodic, its trajectory would contain a cycle $\Omega$,
which would contradict the assumption that $V$ contains no cycles.
In particular, note that $x$'s trajectory contains a cycle if and
only if there is some $y\in Y$ so that $T^{\circ n}\left(x\right)=T^{\circ m}\left(x\right)$
for two distinct positive integers $m$ and $n$; that is, cycles
happen when our trajectories arrive at a point they have already visited.
Thus, writing $x_{n}$ to denote $T^{\circ n}\left(x\right)$, it
must be that $x_{n}\neq x_{m}$ for all distinct non-negative integers
$m,n$.

So, by the pigeonhole principle, given any $R>0$, there must exist
an integer $n_{R}\geq0$ so that $\left|x_{n}\right|_{X}>R$ for all
$n\geq n_{R}$. If not, then there would be infinitely many elements
of $x$'s trajectory which lay in the set $\left\{ y\in Y:\left|y\right|_{X}\leq R\right\} $.
Since $Y$ is discrete, there are only finitely many elements in this
set, which means that the trajectory of $x$ must retrace its steps
at some point: there are distinct $m,n$ for which $x_{m}=x_{n}$,
which, as we saw, is impossible. Since $R$ was arbitrary, we conclude
that $\lim_{n\rightarrow\infty}\left|x_{n}\right|_{X}=\infty$, which
proves that $x$ is a divergent point, and makes $V$ into a divergent
trajectory of $T$. Since $x$ was arbitrary, we conclude that every
element of $V$ is a divergent point of $T$.

Q.E.D.

\section{\label{sec:Hydra-Maps}Hydra Maps}

It is a standard result of algebraic number theory that one can choose
a finite set $\left\{ b_{1},\ldots,b_{N}\right\} $ in $\mathcal{O}_{K}$
so that every element of $\mathcal{O}_{K}$ can be uniquely written
as: 
\begin{equation}
m_{1}b_{1}+\ldots+m_{N}b_{N}
\end{equation}
for some $m_{1},\ldots,m_{N}\in\mathbb{Z}$. That is, as an abelian
group under addition, $\mathcal{O}_{K}$ is a \textbf{lattice}. In
a similar way, every non-zero ideal of $\mathcal{O}_{K}$ can be realized
as lattice, as well.

We use the term \textbf{affine lattice }to refer to a coset of a lattice
in $\mathcal{O}_{K}$; thus, every affine lattice in $\mathcal{O}_{K}$
is of the form: 
\begin{equation}
c+\Lambda\overset{\textrm{def}}{=}\left\{ c+z:z\in\Lambda\right\} 
\end{equation}
for some ideal $\Lambda\subseteq\mathcal{O}_{K}$. We call $\Lambda$
the \textbf{underlying ideal }of $c+\Lambda$. Note that the underlying
ideal is necessarily unique, and is given by $\Lambda=\left\{ z-w:z,w\in c+\Lambda\right\} $.

In this context, letting $\Lambda$ be a lattice in $\mathcal{O}_{K}$
with $p$ distinct cosets, a $p$-Hydra map $H:\mathcal{O}_{K}\rightarrow\mathcal{O}_{K}$
is given by a collection of $p$ affine linear maps, called \textbf{branches},
which we associate to each of the cosets of $\Lambda$. In particular,
given $z\in\mathcal{O}_{K}$, $H\left(z\right)$ is obtained by applying
to $z$ the branch specified by the coset of $\Lambda$ that $z$
happens to lie in. 
\begin{example}
The lattice in question is the ideal $2\mathbb{Z}$, and the parity
of $n$ is precisely the coset of $2\mathbb{Z}$ that $n$ happens
to lie in. 
\end{example}
One difficulty that arises in generalizing maps of this shape to $\mathcal{O}_{K}$
is that while every map of the form: 
\begin{equation}
z\in K\mapsto az+b\in K
\end{equation}
for fixed constants $a,b\in K$ is affine linear, not every affine
linear map on $K$ will be of this form, as there exist linear maps
$K\rightarrow K$ which are not given by multiplication against some
fixed scalar. To deal with this, we introduce the notion of an \textbf{affine
lattice morphism}. 
\begin{defn}
We write $\textrm{Aff}_{\mathbb{F}}K$ to denote the \textbf{group
of all} \textbf{invertible affine $\mathbb{F}$-linear maps} $K\rightarrow K$;
these are maps of the form: 
\begin{equation}
z\mapsto rz+c\label{eq:affine map}
\end{equation}
where $r\in\textrm{End}_{\mathbb{F}}^{\times}\left(K\right)$ and
$c\in K$; $\textrm{Aff}_{\mathbb{F}}K$ is a non-abelian group under
map composition. We write $\left(r,c\right)$ as a shorthand for a
map of the form (\ref{eq:affine map}). 
\end{defn}
\begin{rem}
It is important to distinguish between $\mathbb{Q}$-linear maps on
$K$ and $K$-linear maps on $K$. For the case $K=\mathbb{Q}\left(\sqrt{2}\right)$,
letting $a+b\sqrt{2}\in K$ be arbitrary (where $a,b\in\mathbb{Q}$),
consider the maps $m_{1},m_{2}:K\rightarrow K$ given by:
\begin{align}
m_{1}\left(a+b\sqrt{2}\right)\overset{\textrm{def}}{=} & \sqrt{2}\left(a+b\sqrt{2}\right)\\
m_{2}\left(a+b\sqrt{2}\right)\overset{\textrm{def}}{=} & a+b
\end{align}
Both $m_{1}$ and $m_{2}$ are $\mathbb{Q}$-linear, but only $m_{1}$
is $K$-linear. Indeed:
\begin{align}
m_{2}\left(\sqrt{2}\left(a+b\sqrt{2}\right)\right) & =m_{2}\left(2b+a\sqrt{2}\right)=a+2b\\
\sqrt{2}m_{2}\left(a+b\sqrt{2}\right) & =\sqrt{2}\left(a+b\right)\neq a+2b
\end{align}
\end{rem}
Now we can give the definition of a Hydra map. We begin with the definition
of a \textbf{pre-Hydra}.
\begin{defn}
A pre-Hydra consists of the following data:

I. A global field $K$.

II. A proper, non-zero ideal $\Lambda\subset\mathcal{O}_{K}$.

II. For each element $j\in\mathcal{O}_{K}/\Lambda$, an affine linear
map $\left(r_{j},c_{j}\right)\in\textrm{Aff}_{\mathbb{F}}K$. We often
denote $\left(r_{j},c_{j}\right)$ as $H_{j}$, and call $H_{j}$
the \textbf{$j$th branch} of the pre-Hydra.
\end{defn}
\begin{rem}
By linear algebra, we can always identify $K$ with $\mathbb{F}^{d}$
by choosing an $\mathbb{F}$-basis of $K$ and writing $K$'s elements
in coordinates. Because of this, everything done here can also be
formulated over $\mathbb{F}^{d}$. In that particular case, absolute
values on $K$ then become norms on $\mathbb{F}^{d}$; places of $K$
then become equivalence classes of norms on $\mathbb{F}^{d}$. 
\end{rem}
\begin{rem}
Though we only consider affine linear maps, the formalism extends
just as well to polynomial maps, rational maps, and even analytic
maps. I intend to include these in the formulation in a later version
of this paper. 
\end{rem}
\begin{defn}
A \textbf{$\Lambda$-Hydra (map)} on $K$ is a map $H:\mathcal{O}_{K}\rightarrow\mathcal{O}_{K}$
defined by: 
\begin{equation}
H\left(z\right)\overset{\textrm{def}}{=}\sum_{j\in\mathcal{O}_{K}/\Lambda}\left[z\overset{\Lambda}{\equiv}j\right]\left(r_{j}z+c_{j}\right)\label{eq:Hydra}
\end{equation}
Here, $\left[z\overset{\Lambda}{\equiv}j\right]$ is the function
which is $1$ if $z$ is congruent to $j$ mod $\Lambda$, and is
$0$ otherwise, and $\left(K,\Lambda,\left\{ \left(r_{j},c_{j}\right)\right\} _{j\in\mathcal{O}_{K}/\Lambda}\right)$
is a pre-Hydra so that, for all $j\in\mathcal{O}_{K}/\Lambda$, we
have $r_{j}z+c_{j}\in\mathcal{O}_{K}$ for all $z\in j+\Lambda$. 

We call $\Lambda$ the \textbf{underlying ideal }of $H$, and write
$H_{j}$ to denote the affine linear map $\left(r_{j},c_{j}\right)$,
which we call the \textbf{$j$th branch of $H$}.
\end{defn}
\begin{rem}
As defined, $H\left(z\right)$ applies $r_{j}z+c_{j}$, where $j$
is the representative of the coset of $\Lambda$ to which $z$ belongs.
(\ref{eq:Hydra}) can be written in a single line as: 
\begin{equation}
H\left(z\right)=r_{\left[z\right]_{\Lambda}}z+c_{\left[z\right]_{\Lambda}}
\end{equation}
where $\left[z\right]_{\Lambda}$ denotes the image of $z$ under
the canonical projection $\mathcal{O}_{K}\rightarrow\mathcal{O}_{K}/\Lambda$. 
\end{rem}
\begin{rem}
As $\Lambda$ always has at least two cosets in $\mathcal{O}_{K}$,
we say a map $H:\mathcal{O}_{K}\rightarrow\mathcal{O}_{K}$ is a \textbf{$p$-Hydra
map} whenever $H$ is a $\Lambda$-Hydra map for some $\Lambda\subset\mathcal{O}_{K}$
so that $\left|\mathcal{O}_{K}/\Lambda\right|=p$.
\end{rem}
Some notable properties of $\Lambda$-Hydras are as follows:
\begin{defn}
Let $H$ be a $\Lambda$-Hydra on $K$. We say $H$ is:

I. \textbf{Integral}, if for all $z\in\mathcal{O}_{K}$ and all $j\in\mathcal{O}_{K}/\Lambda$,
$H_{j}\left(z\right)\in\mathcal{O}_{K}$ \emph{if and only if} $z\in j+\Lambda$.

II. \textbf{Proper}, if the map $1-r_{0}\in\textrm{End}_{\mathbb{F}}\left(K\right)$
is invertible.

III. \textbf{Centered}, if $H_{0}\left(0\right)=0$.
\end{defn}
\begin{rem}
From the perspective of dynamical systems, the integrality property
of a Hydra is, by far, its most important property, and plays a critical
role in the proof of the \textbf{Correspondence Principle }(\textbf{Theorem
\ref{thm:CP}}, given at the end of\textbf{ Section \ref{sec:The-Numen-of}}).

The classic Collatz map $C:\mathbb{Z}\rightarrow\mathbb{Z}$ defined
by:
\[
C\left(n\right)\overset{\textrm{def}}{=}\begin{cases}
\frac{n}{2} & \textrm{if }n\overset{2}{\equiv}0\\
3n+1 & \textrm{if }n\overset{2}{\equiv}1
\end{cases}
\]
is a $2$-Hydra on $\mathbb{Z}$ which is non-integral, as the ``odd''
branch $H_{1}\left(x\right)=3x+1$ sends both even and odd integers
to integers. However, this is not a significant obstruction; the classic
``shortening'' of $C$ given by the aptly-named \textbf{Shortened
Collatz Map} $T_{3}:\mathbb{Z}\rightarrow\mathbb{Z}$:
\begin{equation}
T_{3}\left(n\right)\overset{\textrm{def}}{=}\begin{cases}
\frac{n}{2} & \textrm{if }n\textrm{ is even}\\
\frac{3n+1}{2} & \textrm{if }n\textrm{ is odd}
\end{cases}
\end{equation}
\emph{is }integral, as its ``even branch'' $H_{0}$ is integer-valued
only for even integers, while its ``odd'' branch $H_{1}$ is integer-valued
only for odd integers.
\end{rem}
\begin{rem}
Letting $a\in\mathcal{O}_{K}$ be undetermined, one can often conjugate
an uncentered Hydra $H$ to a centered Hydra $H_{\textrm{center}}$
by writing:
\[
H_{\textrm{center}}\left(z\right)=H\left(z+a\right)-a
\]
for a judicious choice of $a\in\mathcal{O}_{K}$.
\end{rem}
\begin{rem}
The centeredness plays a role in the \textbf{Correspondence Principle}.
\end{rem}

\section{\label{sec:The-Numen-of}The Numen of a Hydra Map}

Throughout this section, let $H:\mathcal{O}_{K}\rightarrow\mathcal{O}_{K}$
be a $\Lambda$-Hydra map, not necessarily proper, integral, or centered.
Following \cite{my dissertation,first blog paper}, we construct $X_{H}$,
the numen of $H$, like so. 
\begin{defn}
Let $\Sigma_{\Lambda}^{*}$ denote the set of \textbf{strings} (finite
sequences) of elements of $\mathcal{O}_{K}/\Lambda$ including the
empty string, represented by the empty set $\varnothing$. given any
string $\mathbf{j}=\left(j_{1},\ldots,j_{N}\right)\in\Sigma_{\Lambda}^{*}$,
we call $N$ the \textbf{length }of $\mathbf{j}$, and denote it by
$\left|\mathbf{j}\right|$. The empty string is defined as the unique
string of length $0$. We write $\wedge:\Sigma_{\Lambda}^{*}\times\Sigma_{\Lambda}^{*}\rightarrow\Sigma_{\Lambda}^{*}$
to denote the concatenation operator:
\begin{equation}
\mathbf{i}\wedge\mathbf{j}\overset{\textrm{def}}{=}\left(i_{1},\ldots,i_{\left|\mathbf{i}\right|},j_{1},\ldots,j_{\left|\mathbf{j}\right|}\right)
\end{equation}
We write $\Sigma_{\Lambda}^{\infty}$ to denote the set of strings
of finite or infinite length. Note that we can extend $\wedge$ to
a map $\wedge:\Sigma_{\Lambda}^{*}\times\Sigma_{\Lambda}^{\infty}\rightarrow\Sigma_{\Lambda}^{\infty}$. 
\end{defn}
\begin{defn}
Given any string $\mathbf{j}\in\Sigma_{\Lambda}^{*}$, we define the
\textbf{composition sequence }$H_{\mathbf{j}}\in\textrm{Aff}K$ by:
\begin{equation}
H_{\mathbf{j}}\left(z\right)\overset{\textrm{def}}{=}\left(H_{j_{1}}\circ\cdots\circ H_{j_{\left|\mathbf{j}\right|}}\right)\left(z\right),\textrm{ }\forall z\in K\label{eq:stringmap}
\end{equation}
When $\mathbf{j}$ is the empty string, we define $H_{\varnothing}$
to be the identity map on $K$. 

As an element of $\textrm{Aff}K$, the map $H_{\mathbf{j}}$ can be
written in the form $\left(r,c\right)$ for some $r\in\textrm{End}_{\mathbb{F}}^{\times}\left(K\right)$
and $c\in K$ depending on $\mathbf{j}$. We denote these by $M_{H}\left(\mathbf{j}\right)$
and $X_{H}\left(\mathbf{j}\right)$, respectively, so that: 
\begin{equation}
H_{\mathbf{j}}\left(z\right)=M_{H}\left(\mathbf{j}\right)z+X_{H}\left(\mathbf{j}\right),\textrm{ }\forall z\in K
\end{equation}
In this way, we can realize $M_{H}$ and $X_{H}$ as maps $\Sigma_{\Lambda}^{*}\rightarrow\textrm{End}_{\mathbb{F}}^{\times}\left(K\right)$
and $\Sigma_{\Lambda}^{*}\rightarrow K$, respectively. $X_{H}$ is
called the \textbf{numen of $H$}, and is defined by: 
\begin{equation}
X_{H}\left(\mathbf{j}\right)\overset{\textrm{def}}{=}H_{\mathbf{j}}\left(0\right),\textrm{ }\forall\mathbf{j}\in\Sigma_{\Lambda}^{*}
\end{equation}
\end{defn}
The next two identities are trivial, and follow immediately from the
definitions:
\begin{prop}
For any $\mathbf{i},\mathbf{j}\in\Sigma_{\Lambda}^{*}$:
\begin{equation}
M_{H}\left(\mathbf{i}\wedge\mathbf{j}\right)=M_{H}\left(\mathbf{i}\right)M_{H}\left(\mathbf{j}\right)
\end{equation}
\begin{equation}
X_{H}\left(\mathbf{i}\wedge\mathbf{j}\right)=H_{\mathbf{i}}\left(X_{H}\left(\mathbf{j}\right)\right)
\end{equation}
\end{prop}
\begin{defn}
Letting $p$ denote $\left|\mathcal{O}_{K}/\Lambda\right|$, we fix
\emph{once and for all }choice of an enumeration $\Lambda_{0},\ldots,\Lambda_{p-1}$
of the cosets of $\Lambda$ in $\mathcal{O}_{K}$, with $\Lambda_{0}$
denoting $\Lambda$ itself and $\Lambda_{1},\ldots,\Lambda_{p-1}$
denoting its non-trivial cosets. We then let $\beta:\mathcal{O}_{K}/\Lambda\rightarrow\left\{ 0,\ldots,p-1\right\} $
be the map which sends $i\in\mathcal{O}_{K}/\Lambda$ to the unique
integer $\beta\left(i\right)\in\left\{ 0,\ldots,p-1\right\} $ so
that $i\in\Lambda_{\beta\left(i\right)}$. (Thus, if $i\in\Lambda_{0}$,
then $\beta\left(i\right)=0$; if $i\in\Lambda_{1}$, then $\beta\left(i\right)=1$,
etc.).

Because of this, we shall write $\Sigma_{p}^{*}$ instead of $\Sigma_{\Lambda}^{*}$,
with $\Sigma_{p}^{*}$ denoting the set of all finite strings of the
numbers $\left\{ 0,\ldots,p-1\right\} $, including the empty string.
$\beta$ then induces a bijection between $\Sigma_{\Lambda}^{*}$
and $\Sigma_{p}^{*}$.

Letting $\mathbb{N}_{0}=\left\{ 0,1,2,3,\ldots\right\} $ denote the
\textbf{set of non-negative integers}, we then write $\textrm{DigSum}_{p}:\Sigma_{p}^{*}\rightarrow\mathbb{N}_{0}$
to denote the function: 
\begin{equation}
\textrm{DigSum}_{p}\left(\mathbf{j}\right)=\sum_{n=1}^{\left|\mathbf{j}\right|}\beta\left(j_{n}\right)p^{n-1}\label{eq:DigSum}
\end{equation}
where $\left|\mathbf{j}\right|$ is the \textbf{length }of $\mathbf{j}$
(the number of entries of $\mathbf{j}$). 
\end{defn}
\begin{rem}
Like with writing $\Sigma_{p}^{*}$ instead of $\Sigma_{\Lambda}^{*}$,
going forward, in an very important\textbf{ }abuse of notation, we
will index indexing $r_{j}$ and $c_{j}$ using $j\in\left\{ 0,\ldots,p-1\right\} $
by way of $\beta$, instead of using $j\in\mathcal{O}_{K}/\Lambda$.
Thus, for example, by $r_{1}$, we actually mean $r_{\beta^{-1}\left(1\right)}$,
where $\beta^{-1}\left(1\right)$ is the unique $j\in\mathcal{O}_{K}/\Lambda$
that $\beta$ sends to $1$. 
\end{rem}
\begin{example}
If $\mathcal{O}_{K}=\mathbb{Z}$ and $\Lambda=3\mathbb{Z}$, then
upon identifying $\mathcal{O}_{K}/\Lambda=\mathbb{Z}/3\mathbb{Z}$
with $\left\{ 0,1,2\right\} $, we have that: 
\begin{equation}
\textrm{DigSum}_{3}\left(\mathbf{j}\right)=\sum_{n=1}^{\left|\mathbf{j}\right|}j_{n}3^{n-1}
\end{equation}
Thus, for example: 
\begin{equation}
\textrm{DigSum}_{3}\left(\left(1,2,0,1\right)\right)=1+2\cdot3+0\cdot3^{2}+1\cdot3^{3}
\end{equation}
\end{example}
Although $\textrm{DigSum}_{p}:\Sigma_{p}^{*}\rightarrow\mathbb{N}_{0}$
is surjective, note that it is \emph{not }injective. Indeed, given
any $\mathbf{i}\in\Sigma_{p}^{*}$ if we append finitely many $0$s
to the right of $\mathbf{i}$, we obtain a string $\mathbf{j}\neq\mathbf{i}$
for which $\textrm{DigSum}_{p}\left(\mathbf{i}\right)=\textrm{DigSum}_{p}\left(\mathbf{j}\right)$.
To that end: 
\begin{defn}
Let $\sim$ be the equivalence relation on $\Sigma_{p}^{*}$ defined
by $\mathbf{i}\sim\mathbf{j}$ if and only if $\textrm{DigSum}_{p}\left(\mathbf{i}\right)=\textrm{DigSum}_{p}\left(\mathbf{j}\right)$.
We then write $\Sigma_{p}^{*}/\sim$ to denote the set of equivalence
classes in $\Sigma_{p}^{*}$ under $\sim$. 
\end{defn}
This gives the easy result: 
\begin{prop}
\label{prop:~ equivalence relation}The map $\Sigma_{p}^{*}/\sim\rightarrow\mathbb{N}_{0}$
induced by $\textrm{DigSum}_{p}$ is a bijection. 
\end{prop}
With this, we can realize $X_{H}$ as a function $\mathbb{N}_{0}\rightarrow K$. 
\begin{prop}
If $H$ is centered, then $X_{H}\left(\mathbf{i}\right)=X_{H}\left(\mathbf{j}\right)$
for all $\mathbf{i}$ and $\mathbf{j}$ belonging to the same equivalence
class in $\Sigma_{p}^{*}/\sim$. Hence, $X_{H}:\Sigma_{p}^{*}/\sim\rightarrow K$
is well-defined. We can then view $X_{H}$ as a function $\mathbb{N}_{0}\rightarrow K$
by using $\textrm{DigSum}_{p}$ to identify $\Sigma_{p}^{*}/\sim$
with $\mathbb{N}_{0}$. 
\end{prop}
\begin{rem}
If $H$ is \emph{not} centered, one can still formally define $X_{H}$
without reference to strings by using the functional equations given
in \textbf{Proposition \ref{prop:functional equations}}, below, provided
the functional equations admit a solution.
\end{rem}
Using the definition of $X_{H}\left(\mathbf{j}\right)$ as $H_{\mathbf{j}}\left(0\right)$,
one can easily give the following functional equation characterization
of $X_{H}:\mathbb{N}_{0}\rightarrow K$.
\begin{prop}[Siegel (2022); see \cite{my dissertation}]
\label{prop:functional equations}The set of functions $X_{H}:\mathbb{N}_{0}\rightarrow K$
satisfying the functional equations: 
\begin{equation}
X_{H}\left(pn+j\right)=r_{j}X_{H}\left(n\right)+c_{j},\textrm{ }\forall n\in\mathbb{N}_{0},\textrm{ }\forall j\in\left\{ 0,\ldots,p-1\right\} \label{eq:X_H functional equation on N_0}
\end{equation}
(where $r_{j}X_{H}\left(n\right)$ is the image of $X_{H}\left(n\right)\in K$
under $r_{j}\in\textrm{End}_{\mathbb{F}}^{\times}\left(K\right)$)
is in a bijective correspondence with the set of $z\in K$ for which
$\left(1-r_{0}\right)z=c_{0}$. In particular:

I. If $H$ is proper, then (\ref{eq:X_H functional equation on N_0})
has a unique solution, satisfying $X_{H}\left(0\right)=\left(1-r_{0}\right)^{-1}c_{0}$.

II. If $H$ is non-proper, then (\ref{eq:X_H functional equation on N_0})
has a solution $X_{H}$ if and only if there are $z\in K$ for which
$\left(1-r_{0}\right)z=c_{0}$. For each such $z$, one can set $X_{H}\left(0\right)=z$
and then obtain a distinct solution of (\ref{eq:X_H functional equation on N_0}).
In this way, we make think of $X_{H}\left(0\right)=z$ as an initial
condition subject to which (\ref{eq:X_H functional equation on N_0})
can be solved.
\end{prop}
Next, we introduce the ring of $p$-adic integers $\mathbb{Z}_{p}$.
\begin{defn}
Given an integer $p\geq2$, the ring of $p$-adic integers is defined
by the following quotient of a formal power series ring: 
\begin{equation}
\mathbb{Z}_{p}\overset{\textrm{def}}{=}\frac{\mathbb{Z}\left[\left[x\right]\right]}{\left\langle x-p\right\rangle }
\end{equation}
In practice, one can think of a $p$-adic integer as an integer written
in base $p$, but with the caveat that said integer can potentially
have infinitely many base $p$ digits. The quotient construction justifies
viewing $p$-adic integers as ``formal power series in $p$ with
coefficients in $\left\{ 0,\ldots,p-1\right\} $''. Given $\mathfrak{z}\in\mathbb{Z}_{p}$,
this power series representation of $\mathfrak{z}$ is called the
\textbf{Hensel series of $\mathfrak{z}$}, and the coefficients of
this series are called the \textbf{$p$-adic digits} \textbf{of $\mathfrak{z}$}.
$\mathbb{Z}_{p}$ is made into a compact abelian group by equipping
it with the adic topology induced by the maximal ideal $x-p$ of $\mathbb{Z}\left[\left[x\right]\right]$.

When $p$ is a prime number, $\mathbb{Z}_{p}$ can be topologized
using the \textbf{$p$-adic absolute value}. This is obtained by defining
the \textbf{$p$-adic valuation }$v_{p}:\mathbb{Q}\rightarrow\mathbb{Z}\cup\left\{ +\infty\right\} $
by: 
\begin{equation}
v_{p}\left(n\right)=\#\textrm{ of times }p\textrm{ divides }n
\end{equation}
for all positive integers $n$, and then using the identities: 
\begin{equation}
v_{p}\left(\frac{1}{n}\right)=-v_{p}\left(n\right),\textrm{ }\forall n\geq1
\end{equation}
and 
\begin{equation}
v_{p}\left(rs\right)=v_{p}\left(r\right)+v_{p}\left(s\right),\textrm{ }\forall r,s\in\mathbb{Q}
\end{equation}
to extend $v_{p}$ to $\mathbb{Q}$. Note that $v_{p}\left(0\right)$
is defined to be $+\infty$. The $p$-adic absolute value is then
defined by: 
\begin{equation}
\left|r\right|_{p}\overset{\textrm{def}}{=}p^{-v_{p}\left(r\right)},\textrm{ }\forall r\in\mathbb{Q}
\end{equation}
In this way, $\mathbb{Z}_{p}$ can be constructed as the completion
of the metric space obtained by equipping $\mathbb{Z}$ with the metric
induced by the $p$-adic absolute value. $\mathbb{Q}_{p}$, the field
of $p$-adic numbers, is the field of fractions of $\mathbb{Z}_{p}$,
or, equivalently, the completion of $\mathbb{Q}$ with respect to
the metric induced by $\left|\cdot\right|_{p}$. Importantly, the
$p$-adic absolute value is \emph{non-archimedean}, meaning it satisfies
the \textbf{ultrametric inequality}: 
\begin{equation}
\left|\mathfrak{x}+\mathfrak{y}\right|_{p}\leq\max\left\{ \left|\mathfrak{x}\right|_{p},\left|\mathfrak{y}\right|_{p}\right\} ,\textrm{ }\forall\mathfrak{x},\mathfrak{y}\in\mathbb{Q}_{p}
\end{equation}
with equality whenever $\left|\mathfrak{x}\right|_{p}\neq\left|\mathfrak{y}\right|_{p}$.

Finally, given any $\mathfrak{z}\in\mathbb{Z}_{p}$ and any integer
$n\geq0$, we write $\left[\mathfrak{z}\right]_{p^{n}}$ to denote
the unique integer in $\left\{ 0,\ldots,p^{n}-1\right\} $ which is
congruent to $\mathfrak{z}$ mod $p^{n}$. 
\end{defn}
\begin{rem}
Using this definition, it can be shown that, as a compact abelian
group, $\mathbb{Z}_{p}$ is isomorphic to the direct product: 
\begin{equation}
\prod_{\ell\mid p}\mathbb{Z}_{\ell}
\end{equation}
taken over all prime divisors $\ell$ of $p$. In particular, this
shows that for any prime $\ell$ and any integer $n\geq1$, $\mathbb{Z}_{\ell^{n}}$
is isomorphic to $\mathbb{Z}_{\ell}$ as a topological ring. As such,
we will write the above direct product isomorphism as: 
\begin{equation}
\mathbb{Z}_{p}\cong\prod_{\ell\mid p}\mathbb{Z}_{\ell^{v_{\ell}\left(p\right)}}
\end{equation}
Ex: 
\begin{equation}
\mathbb{Z}_{12}\cong\mathbb{Z}_{3}\times\mathbb{Z}_{4}
\end{equation}
\end{rem}
Just like with $\mathbb{N}_{0}$ and $\Sigma_{\Lambda}^{*}$, we can
identify $\mathbb{Z}_{p}$ with a set of strings.
\begin{defn}
We write:

I. $\Sigma_{p}^{\infty}$, to denote the set of all strings of elements
of\emph{ }$\left\{ 0,\ldots,p-1\right\} $ of either finite or infinite
length. We also include the empty string as an element of $\Sigma_{p}^{\infty}$.

II. $\Sigma_{p}^{+\infty}$, to denote the set of all strings of elements
of\emph{ }$\left\{ 0,\ldots,p-1\right\} $ of infinite length. 
\end{defn}
\begin{rem}
We can and will extend $\textrm{DigSum}_{p}$ from a function $\Sigma_{p}^{*}\rightarrow\mathbb{N}_{0}$
to a function $\Sigma_{p}^{\infty}\rightarrow\mathbb{Z}_{p}$, where
$\Sigma_{p}^{\infty}$ is the set of finite- or infinite-length strings
in the symbols $\left\{ 0,\ldots,p-1\right\} $. Letting $\sim$ denote
the equivalence relation on $\Sigma_{p}^{\infty}$ with $\mathbf{i}\sim\mathbf{j}$
if and only if $\textrm{DigSum}_{p}\left(\mathbf{i}\right)=\textrm{DigSum}_{p}\left(\mathbf{j}\right)$,
we then have that $\textrm{DigSum}_{p}$ induces a bijection $\Sigma_{p}^{\infty}/\sim\rightarrow\mathbb{Z}_{p}$,
where $\Sigma_{p}^{\infty}/\sim$ is the set of equivalence classes
of $\Sigma_{p}^{\infty}$ under $\sim$. 
\end{rem}
In \cite{my dissertation}, it was shown that, given mild assumptions
on the $r_{j}$s, one can use the identification of $\Sigma_{p}^{\infty}/\sim$
with $\mathbb{Z}_{p}$ to interpolate $X_{H}:\mathbb{N}_{0}\rightarrow K$
to a function on $\mathbb{Z}_{p}$. The details for this are given
in \textbf{Lemma \ref{lem:extension and measurability of X_H}}, below.
The essential idea is that by making a judicious choice of absolute
values on $K$, one can realize $X_{H}$ as a function (or, in general,
as a distribution) out of $\mathbb{Z}_{p}$ whose regularity depends
on how close the $r_{j}$s are to behaving as contractions with respect
to $K$. The behavior of the $r_{j}$s is captured by their norm with
respect to a given absolute value of $K$. The result given below
is a slightly modified version of the versions originally stated in
\cite{first blog paper}. The reader should be aware that \textbf{Lemma
\ref{lem:extension and measurability of X_H} }treats $\mathbb{Z}_{p}$
as a measure space equipped with its \textbf{real-valued Haar probability
measure}, and treat $K_{\ell}$ as a measure space equipped with its
real-valued Haar measure. We also use the concept of a \textbf{locally
constant function}. Readers unfamiliar with either of these can find
an expository treatment given of them at the start of \textbf{Section
\ref{sec:-Adic-Fourier-Analysis}}.
\begin{lem}
\label{lem:extension and measurability of X_H}Suppose $H$ is proper.
For each non-trivial place $\ell$ of $K$, let:

\begin{equation}
\rho_{H,\ell}\overset{\textrm{def}}{=}\prod_{j=0}^{p-1}\left\Vert r_{j}\right\Vert _{\ell}
\end{equation}
be the product of the $\ell$-adic norms of the $r_{j}$s. 

Now, consider the extension $X_{H}:\mathbb{Z}_{p}\rightarrow K_{\ell}$
defined by the formula:
\begin{align}
X_{H}\left(\mathfrak{z}\right) & \overset{K_{\ell}}{=}\lim_{n\rightarrow\infty}X_{H}\left(\left[\mathfrak{z}\right]_{p^{n}}\right)\label{eq:extension of Chi_H}
\end{align}
for all $\mathfrak{z}\in\mathbb{Z}_{p}$ for which $X_{H}\left(\mathfrak{z}\right)$
exists. Then:

I. If $\rho_{H,\ell}<1$, then (\ref{eq:extension of Chi_H}) defines
a measurable function $X_{H}:\mathbb{Z}_{p}\rightarrow K_{\ell}$,
with (\ref{eq:extension of Chi_H}) converging for almost every $\mathfrak{z}\in\mathbb{Z}_{p}$.
In particular, it converges at all $\mathfrak{z}\in\mathbb{Z}_{p}$
for which:
\begin{equation}
\sum_{j=0}^{p-1}d_{H,\ell,j}\left(\mathfrak{z}\right)\ln\left\Vert r_{j}\right\Vert _{\ell}<0\label{eq:convergence condition}
\end{equation}
where: 
\begin{equation}
d_{H,\ell,j}\left(\mathfrak{z}\right)\overset{\textrm{def}}{=}\begin{cases}
\limsup_{m\rightarrow\infty}\frac{\#_{j}\left(\left[\mathfrak{z}\right]_{p^{m}}\right)}{m} & \textrm{if }\left\Vert r_{j}\right\Vert _{\ell}\leq1\\
\liminf_{m\rightarrow\infty}\frac{\#_{j}\left(\left[\mathfrak{z}\right]_{p^{m}}\right)}{m} & \textrm{if }\left\Vert r_{j}\right\Vert _{\ell}>1
\end{cases}\label{eq:upper density of j}
\end{equation}
is the \textbf{$\left(H,\ell\right)$-normalized upper density of
$j$ in the $p$-adic digits of $\mathfrak{z}$}.

II. If $\rho_{H,\ell}<1$ and $\max_{0\leq j<p}\left\Vert r_{j}\right\Vert _{\ell}<1$,
then (\ref{eq:extension of Chi_H}) converges uniformly with respect
to $\mathfrak{z}\in\mathbb{Z}_{p}$, making $X_{H}:\mathbb{Z}_{p}\rightarrow K_{\ell}$
into a continuous function.

Moreover, given either (I) or (II), we have that:
\begin{equation}
X_{H}\left(p\mathfrak{z}+j\right)=r_{j}X_{H}\left(\mathfrak{z}\right)+c_{j}\label{eq:X_H functional equation}
\end{equation}
occurs for all $j\in\left\{ 0,\ldots,p-1\right\} $ and all $\mathfrak{z}\in\mathbb{Z}_{p}$
for which $X_{H}\left(\mathfrak{z}\right)$ converges in $K_{\ell}$,
and that the extensions of $X_{H}$ given in (I) and (II), respectively,
are the unique functions satisfying these functional equations. Furthermore,
if $\ell$ is non-archimedean and $\max_{0\leq j<p}\left\Vert r_{j}\right\Vert _{\ell}\leq1$,
then:
\begin{equation}
\limsup_{N\rightarrow\infty}\sup_{\mathfrak{z}\in\mathbb{Z}_{p}}\left|X_{H}\left(\left[\mathfrak{z}\right]_{p^{N}}\right)\right|_{\ell}\leq\max_{0\leq j<p}\left|c_{j}\right|_{\ell}\label{eq:a.e.  ell-adic bound on X_H in NA case}
\end{equation}
In particular, $X_{H}:\mathbb{Z}_{p}\rightarrow K_{\ell}$ is $\ell$-adically
bounded by $\max_{0\leq j<p}\left|c_{j}\right|_{\ell}$ on its set
of convergence.
\end{lem}
\begin{rem}
Like with \textbf{Proposition \ref{prop:functional equations}}, if
$H$ is not proper but there exist $z\in K$ for which $\left(1-r_{0}\right)z=c_{0}$,
one can recover the result of \textbf{Lemma \ref{lem:extension and measurability of X_H}
}subject to the condition that the uniqueness of $X_{H}$ is requires
not only the condition (\ref{eq:extension of Chi_H}) and the functional
equations (\ref{eq:X_H functional equation}), but also that $X_{H}$
satisfy the initial condition $X_{H}\left(0\right)=z$.
\end{rem}
Proof: First, note that if $\mathfrak{z}\in\mathbb{N}_{0}$, then
$\left[\mathfrak{z}\right]_{p^{n}}=\mathfrak{z}$ for all sufficiently
large $n$, and hence (\ref{eq:extension of Chi_H}) converges. So,
suppose $\mathfrak{z}\in\mathbb{Z}_{p}\backslash\mathbb{N}_{0}$.
Next, following \cite{first blog paper}, a simple inductive argument
gives the following formula for the value of $X_{H}$ at a given string
$\mathbf{j}$ of finite length: 
\begin{equation}
X_{H}\left(\mathbf{j}\right)=\sum_{m=0}^{\left|\mathbf{j}\right|-1}r_{\left[\mathbf{j}\right]_{m}}c_{j_{m+1}},\textrm{ }\forall\mathbf{j}\in\Sigma_{p}^{*}\label{eq:Chi H of bold j}
\end{equation}
where: 
\begin{equation}
r_{\left[\mathbf{j}\right]_{m}}\overset{\textrm{def}}{=}\prod_{k=1}^{m}r_{j_{k}}=\begin{cases}
\textrm{Id} & \textrm{if }\mathbf{j}=\varnothing\\
r_{j_{1}}\cdots r_{j_{m}} & \textrm{if }\left|\mathbf{j}\right|\geq1
\end{cases}
\end{equation}

Applying norms, and using 
\begin{equation}
\left|r_{j}z\right|_{\ell}\leq\left\Vert r_{j}\right\Vert _{\ell}\left|z\right|_{\ell},\textrm{ }\forall z\in K,\textrm{ }\forall j\in\left\{ 0,\ldots,p-1\right\} 
\end{equation}
the triangle inequality gives: 
\begin{equation}
\left|X_{H}\left(\mathbf{j}\right)\right|_{\ell}=\sum_{m=0}^{\left|\mathbf{j}\right|-1}\left|c_{j_{m+1}}\right|_{\ell}\prod_{k=1}^{m}\left\Vert r_{j_{k}}\right\Vert _{\ell}
\end{equation}
So, given $\mathfrak{z}\in\mathbb{Z}_{p}$, for each $n\geq0$, let:
\begin{equation}
\mathbf{j}_{n}\left(\mathfrak{z}\right)=\left(j_{1}\left(\mathfrak{z}\right),j_{2}\left(\mathfrak{z}\right),\ldots,j_{n}\left(\mathfrak{z}\right)\right)
\end{equation}
denote the string of $p$-adic digits of $\left[\mathfrak{z}\right]_{p^{n}}$;
that is: 
\begin{equation}
\left[\mathfrak{z}\right]_{p^{n}}=\textrm{DigSum}_{p}\left(\mathbf{j}_{n}\left(\mathfrak{z}\right)\right)
\end{equation}
Thus, we have: 
\begin{equation}
\left|X_{H}\left(\left[\mathfrak{z}\right]_{p^{n}}\right)\right|_{\ell}=\left|X_{H}\left(\mathbf{j}_{n}\left(\mathfrak{z}\right)\right)\right|_{\ell}\ll\sum_{m=0}^{n-1}\prod_{k=1}^{m}\left\Vert r_{j_{k}\left(\mathfrak{z}\right)}\right\Vert _{\ell}=\sum_{m=0}^{n-1}\prod_{j=0}^{p-1}\left\Vert r_{j}\right\Vert _{\ell}^{\#_{j}\left(\left[\mathfrak{z}\right]_{p^{m}}\right)}
\end{equation}
where $\#_{j}\left(\left[\mathfrak{z}\right]_{p^{m}}\right)$ counts
the number of digits of $\left[\mathfrak{z}\right]_{p^{m}}$ which
are equal to $j$. Hence:
\begin{equation}
\left|X_{H}\left(\left[\mathfrak{z}\right]_{p^{n}}\right)\right|_{\ell}\ll\sum_{m=0}^{n-1}\prod_{j=0}^{p-1}\left\Vert r_{j}\right\Vert _{\ell}^{\#_{j}\left(\left[\mathfrak{z}\right]_{p^{m}}\right)}\label{eq:Ready for estimation}
\end{equation}
where the constant of proportionality is $C$.

Applying the \textbf{Root Test} for series convergence, and noting
that $0\leq\rho_{j}\leq1$ for all $j$, we see that the upper bound
in (\ref{eq:Ready for estimation}) will converge provided: 
\begin{equation}
\limsup_{m\rightarrow\infty}\prod_{j=0}^{p-1}\left\Vert r_{j}\right\Vert _{\ell}^{\frac{\#_{j}\left(\left[\mathfrak{z}\right]_{p^{m}}\right)}{m}}<1\label{eq:convergence condition-1}
\end{equation}
Since, for almost every $\mathfrak{z}\in\mathbb{Z}_{p}$ (in the sense
of $\mathbb{Z}_{p}$'s Haar probability measure), one has:
\begin{equation}
\lim_{m\rightarrow\infty}\frac{\#_{j}\left(\left[\mathfrak{z}\right]_{p^{m}}\right)}{m}=\frac{1}{p}
\end{equation}
for all $j\in\left\{ 0,\ldots,p-1\right\} $, (\ref{eq:convergence condition-1})
is equivalent to: 
\begin{equation}
\prod_{j=0}^{p-1}\left\Vert r_{j}\right\Vert _{\ell}^{1/p}<1\label{eq:geometric mean condition}
\end{equation}
The left-hand side is $\rho_{H,\ell}^{1/p}$. Thus, $\rho_{H,\ell}<1$
implies almost-everywhere convergence. As $X_{H}$ is defined by the
limit $\lim_{n\rightarrow\infty}X_{H}\left(\left[\mathfrak{z}\right]_{p^{n}}\right)$,
and since $X_{H}\left(\left[\mathfrak{z}\right]_{p^{n}}\right)$ is
a locally constant, and hence, measurable $K_{\ell}$-valued function
on $\mathbb{Z}_{p}$, it follows that the almost-everywhere convergent
limit:
\begin{equation}
X_{H}\left(\mathfrak{z}\right)\overset{K_{\ell}}{=}\lim_{n\rightarrow\infty}X_{H}\left(\left[\mathfrak{z}\right]_{p^{n}}\right)
\end{equation}
 is measurable, as well. 

Finally, if $\max_{0\leq j<p}\left\Vert r_{j}\right\Vert _{\ell}<1$,
set: 
\begin{equation}
\rho=\max_{0\leq j<p}\left\Vert r_{j}\right\Vert _{\ell}
\end{equation}
Then, we re-write (\ref{eq:Chi H of bold j}) like so:
\begin{equation}
X_{H}\left(\left[\mathfrak{z}\right]_{p^{n}}\right)=\sum_{m=0}^{n-1}\left(\prod_{k=0}^{m-1}r_{\left[\theta_{p}^{\circ k}\left(\mathfrak{z}\right)\right]_{p}}\right)c_{\left[\theta_{p}^{\circ m}\left(\mathfrak{z}\right)\right]_{p}},\textrm{ }\forall\mathfrak{z}\in\mathbb{Z}_{p},\textrm{ }\forall n\geq0\label{eq:Chi_H of z mod p^n}
\end{equation}
Here, $\theta_{p}:\mathbb{Z}_{p}\rightarrow\mathbb{Z}_{p}$ is the
\textbf{$p$-adic shift map}, which acts on Hensel series representations
as:
\begin{equation}
\theta_{p}\left(\sum_{n=0}^{\infty}d_{n}p^{n}\right)\overset{\textrm{def}}{=}\sum_{n=0}^{\infty}d_{n+1}p^{n}
\end{equation}
i.e.:
\begin{equation}
\theta_{p}\left(\mathfrak{z}\right)=\frac{\mathfrak{z}-\left[\mathfrak{z}\right]_{p}}{p}
\end{equation}
$\left[\theta_{p}^{\circ k}\left(\mathfrak{z}\right)\right]_{p}$,
the projection mod $p$ of the image of $\mathfrak{z}$ under $k$
applications of $\theta_{p}$, is then precisely the coefficient of
$p^{k}$ in $\mathfrak{z}$'s Hensel series representation. This then
gives us the formal expression:
\begin{equation}
\lim_{N\rightarrow\infty}X_{H}\left(\left[\mathfrak{z}\right]_{p^{N}}\right)=\sum_{m=0}^{\infty}\left(\prod_{k=0}^{m-1}r_{\left[\theta_{p}^{\circ k}\left(\mathfrak{z}\right)\right]_{p}}\right)c_{\left[\theta_{p}^{\circ m}\left(\mathfrak{z}\right)\right]_{p}}\label{eq:F series for rising continuation of X_H}
\end{equation}
We then have:
\begin{align*}
\left|X_{H}\left(\left[\mathfrak{z}\right]_{p^{n}}\right)-\lim_{N\rightarrow\infty}X_{H}\left(\left[\mathfrak{z}\right]_{p^{N}}\right)\right|_{\ell} & =\left|\sum_{m=n}^{\infty}\left(\prod_{k=0}^{m-1}r_{\left[\theta_{p}^{\circ k}\left(\mathfrak{z}\right)\right]_{p}}\right)c_{\left[\theta_{p}^{\circ m}\left(\mathfrak{z}\right)\right]_{p}}\right|_{\ell}\\
 & \leq\sum_{m=n}^{\infty}\prod_{k=0}^{m-1}\left\Vert r_{\left[\theta_{p}^{\circ k}\left(\mathfrak{z}\right)\right]_{p}}\right\Vert _{\ell}\left|c_{\left[\theta_{p}^{\circ m}\left(\mathfrak{z}\right)\right]_{p}}\right|_{\ell}\\
 & \leq\sum_{m=n}^{\infty}\prod_{k=0}^{m-1}\rho\left(\max_{0\leq j<p}\left|c_{j}\right|_{\ell}\right)\\
 & \ll\sum_{m=n}^{\infty}\rho^{m}\\
\left(\rho<1\right); & \overset{\mathbb{R}}{=}\frac{\rho^{n}}{1-\rho}
\end{align*}
As $n\rightarrow\infty$, the upper bound tends to $0$ uniformly
with respect to $\mathfrak{z}\in\mathbb{Z}_{p}$, thereby establishing
the uniform convergence of $\lim_{N\rightarrow\infty}X_{H}\left(\left[\mathfrak{z}\right]_{p^{N}}\right)$.
As locally constant functions out of $\mathbb{Z}_{p}$ are continuous,
this then proves that:
\begin{equation}
X_{H}\left(\mathfrak{z}\right)\overset{\textrm{def}}{=}\lim_{N\rightarrow\infty}X_{H}\left(\left[\mathfrak{z}\right]_{p^{N}}\right)
\end{equation}
is continuous.

Almost done, observing the identity: 
\begin{equation}
X_{H}\left(\left[p\mathfrak{z}+j\right]_{p^{n}}\right)=X_{H}\left(p\left[\mathfrak{z}\right]_{p^{n-1}}+j\right)=r_{j}X_{H}\left(\left[\mathfrak{z}\right]_{p^{n-1}}\right)+\mathbf{c}_{j}
\end{equation}
for all $n\geq1$, all $\mathfrak{z}\in\mathbb{Z}_{p}$, and all $j\in\left\{ 0,\ldots,p-1\right\} $,
we see that the functional equations (\ref{eq:X_H functional equation})
are satisfied for all $\mathfrak{z}$ for which $\lim_{n\rightarrow\infty}X_{H}\left(\left[\mathfrak{z}\right]_{p^{n}}\right)$
exists.

For the uniqueness, given any function $f:\mathbb{N}_{0}\rightarrow K_{\ell}$
satisfying:
\begin{equation}
f\left(pn+j\right)=r_{j}f\left(n\right)+c_{j},\textrm{ }\forall n\geq0,\textrm{ }\forall j\in\left\{ 0,\ldots,p-1\right\} \label{eq:generic functional equation}
\end{equation}
observe that since $1-r_{0}$ is invertible by the properness of $H$,
we have for $n=j=0$:
\begin{equation}
f\left(0\right)=r_{0}f\left(0\right)+c_{0}
\end{equation}
and hence:
\begin{equation}
f\left(0\right)=\left(1-r_{0}\right)^{-1}c_{0}
\end{equation}
By induction, it follows that this value of $f\left(0\right)$ then
uniquely determines $f\left(n\right)$ for all $n\geq0$. This shows
that the restriction of $X_{H}$ to $\mathbb{N}_{0}$ is unique, and
hence, that the continuation of $X_{H}$ to $\mathbb{Z}_{p}$ is also
unique.

Finally, suppose that $\ell$ is non-archimedean, and that the maximum
of the $\ell$-adic norms of the $r_{j}$s is $\leq1$. Then (\ref{eq:F series for rising continuation of X_H})
and the ultrametric property of the $\ell$-adic norm implies that:
\begin{equation}
\limsup_{N\rightarrow\infty}\left|X_{H}\left(\left[\mathfrak{z}\right]_{p^{N}}\right)\right|_{\ell}\leq\max_{0\leq j<p}\left|c_{j}\right|_{\ell}
\end{equation}
which is (\ref{eq:a.e.  ell-adic bound on X_H in NA case}), and which
establishes the $\ell$-adic boundedness of $X_{H}$.

Q.E.D.

\vphantom{}
\begin{rem}
As constructed, note that for any $\mathfrak{z}\in\mathbb{Z}_{p}\cap\mathbb{Q}$
(that is, for any \emph{rational }$p$-adic integer) for which $X_{H}\left(\mathfrak{z}\right)$
is finite-valued, $X_{H}\left(\mathfrak{z}\right)$ will necessarily
be an element of $K$. Since $X_{H}\left(\mathfrak{z}\right)$ is
the limit of $X_{H}\left(\left[\mathfrak{z}\right]_{p^{n}}\right)$
as $n\rightarrow\infty$, one can use (\ref{eq:Chi H of bold j})
to write:
\begin{equation}
X_{H}\left(\mathfrak{z}\right)=\sum_{m=0}^{\infty}\left(\prod_{k=0}^{m-1}r_{d_{k}\left(\mathfrak{z}\right)}\right)c_{d_{m}\left(\mathfrak{z}\right)}\label{eq:primitive F-series}
\end{equation}
where:
\begin{equation}
\mathfrak{z}=\sum_{n=0}^{\infty}d_{n}\left(\mathfrak{z}\right)p^{n}
\end{equation}
is the Hensel series of $\mathfrak{z}$, and where:
\begin{equation}
\prod_{k=0}^{m-1}r_{d_{k}\left(\mathfrak{z}\right)}\overset{\textrm{def}}{=}r_{d_{0}\left(\mathfrak{z}\right)}\circ\cdots\circ r_{d_{m-1}\left(\mathfrak{z}\right)}
\end{equation}
with the product being the identity map whenever $m=0$. Because $\mathfrak{z}$
is rational if and only if its $p$-adic digits are eventually periodic,
the eventual-periodicity of the $d_{m}$s guarantees that there will
be an $M$, depending on $\mathfrak{z}$, so that the series formed
by summing the terms of (\ref{eq:primitive F-series}) for $m\geq M$
is a $K$-linear combination of finitely many geometric series. If
$X_{H}\left(\mathfrak{z}\right)$ is convergent, these geometric series
will necessarily converge to elements of $K$.
\end{rem}
\begin{rem}
In the study of iterated function systems (IFS), one traditionally
use spaces of strings\textemdash usually called \textbf{\emph{words}}
in that subject\textemdash where, following \cite{my dissertation},
we use the $p$-adic integers. In passing from arbitrary strings $\mathbf{j}\in\Sigma_{p}^{\infty}$
to the $p$-adic integers using $\textrm{DigSum}_{p}$, we lose distinction
between those strings that end in infinitely many $0$s. Given any
iterated function system generated by maps $f_{1},\ldots,f_{N}$ on
a complete normed vector space $V$, we can avoid this lost information
by defining $p$ to be $N+1$, and then setting $H_{0}$ to be the
identity map of $V$ and setting $H_{j}=f_{j}$ for all $j\in\left\{ 1,\ldots,p-1\right\} $.
The numen:
\begin{equation}
X_{H}\left(p\mathfrak{z}+j\right)=H_{j}\left(X_{H}\left(\mathfrak{z}\right)\right),\textrm{ }\forall j\in\left\{ 0,\ldots,p-1\right\} ,\textrm{ }\forall\mathfrak{z}\in\mathbb{Z}_{p}
\end{equation}
is then a surjective map from $\mathbb{Z}_{p}$ onto the fractal attractor
$F$ of the $f_{j}$s. If we set $q=N$ and define: 
\begin{equation}
Y_{H}\left(q\mathfrak{z}+j\right)=f_{j+1}\left(Y_{H}\left(\mathfrak{z}\right)\right),\textrm{ }\forall j\in\left\{ 0,\ldots,q-1\right\} ,\textrm{ }\forall\mathfrak{z}\in\mathbb{Z}_{q}
\end{equation}
the image of $Y_{H}$ would omit all the points of $F$ corresponding
to composition sequences:
\begin{equation}
f_{j_{1}}\circ f_{j_{2}}\circ\cdots
\end{equation}
so that $j_{n}=1$ for all sufficiently large $n$. 
\end{rem}
Finally, we have the analogue of the \textbf{Correspondence Principle}
from \cite{first blog paper}.
\begin{thm}[Correspondence Principle]
\label{thm:CP}Let $H:\mathcal{O}_{K}\rightarrow\mathcal{O}_{K}$
be a proper, centered, integral $\Lambda$-Hydra map. Suppose:

i. For each $n\geq1$, there is a place $\ell_{n}$ of $K$ so that
$\left\Vert M_{H}\left(n\right)\right\Vert _{\ell_{n}}<1$.

ii. There is a non-trivial non-archimedean place $v$ of $K$ so that
$\left\Vert r_{j}\right\Vert _{v}>1$ for all $j\in\left\{ 0,\ldots,p-1\right\} $
and $v\left(\Lambda\right)\geq1$. 

Then, $z\in\mathcal{O}_{K}$ is a periodic point of $H$ if and only
if there is a rational $p$-adic integer $\mathfrak{z}\in\mathbb{Z}_{p}\cap\mathbb{Q}$
so that $X_{H}\left(\mathfrak{z}\right)=z$. 
\end{thm}
Proof (Sketch): Let $\mathbf{j}\in\Sigma_{p}^{*}$ be a string so
that $H_{\mathbf{j}}\left(z\right)=z$, and let $n=\textrm{DigSum}_{p}\left(\mathbf{j}\right)$.
Then: 
\begin{equation}
z=H_{\mathbf{j}}\left(z\right)=M_{H}\left(n\right)z+X_{H}\left(n\right)
\end{equation}
Consequently, for all $m\geq1$:
\begin{equation}
z=H_{\mathbf{j}}^{\circ m}\left(z\right)=M_{q}\left(\mathbf{j}^{\wedge m}\right)z+X_{H}\left(\mathbf{j}^{\wedge m}\right)=\left(M_{q}\left(n\right)\right)^{m}z+X_{H}\left(\mathbf{j}^{\wedge m}\right)
\end{equation}
where $\mathbf{j}^{\wedge m}$ is the string obtained by concatenating
$m$ copies of $\mathbf{j}$. Here:
\begin{equation}
X_{H}\left(\mathbf{j}^{\wedge m}\right)=H_{\mathbf{j}}\left(X_{H}\left(\mathbf{j}^{\wedge m-1}\right)\right)=M_{H}\left(\mathbf{j}\right)X_{H}\left(\mathbf{j}^{\wedge m-1}\right)+X_{H}\left(\mathbf{j}\right)
\end{equation}
from which it follows that:
\begin{equation}
X_{H}\left(\mathbf{j}^{\wedge m}\right)=\frac{1-M_{H}^{m}\left(\mathbf{j}\right)}{1-M_{H}\left(\mathbf{j}\right)}X_{H}\left(\mathbf{j}\right)=\frac{1-M_{H}^{m}\left(n\right)}{1-M_{H}\left(n\right)}X_{H}\left(n\right)
\end{equation}
By (i), we can choose a place $\ell_{n}$ so that $\left\Vert M_{H}^{m}\left(n\right)\right\Vert _{\ell_{n}}=\left\Vert M_{H}\left(n\right)\right\Vert _{\ell_{n}}^{m}$
tends to $0$ as $m\rightarrow\infty$. This gives us:
\begin{equation}
z=H_{\mathbf{j}}^{\circ m}\left(z\right)=\left(M_{q}\left(n\right)\right)^{m}z+\frac{1-M_{H}^{m}\left(n\right)}{1-M_{H}\left(n\right)}X_{H}\left(n\right),\textrm{ }\forall m\geq1
\end{equation}
and hence, we get $\ell$-adic convergence:
\begin{equation}
z=\frac{X_{H}\left(n\right)}{1-M_{H}\left(n\right)}
\end{equation}
as $m\rightarrow\infty$. As in \cite{first blog paper}, we have:
\begin{equation}
\frac{X_{H}\left(n\right)}{1-M_{H}\left(n\right)}=X_{H}\left(B_{p}\left(n\right)\right),\textrm{ }\forall n\geq1
\end{equation}
where:
\begin{equation}
B_{p}\left(n\right)\overset{\textrm{def}}{=}\frac{n}{1-p^{\lambda_{p}\left(n\right)}}
\end{equation}
is the function which outputs the $p$-adic integer whose sequence
of digits is precisely the result of concatenating infinitely many
copies of the sequence of $n$'s $p$-adic digits. Since $B_{p}\left(n\right)\in\left(\mathbb{Z}_{p}\backslash\mathbb{N}_{0}\right)\cap\mathbb{Q}$
for all $n\geq1$, this shows that $z$ is in the image of $\left(\mathbb{Z}_{p}\backslash\mathbb{N}_{0}\right)\cap\mathbb{Q}$
under $X_{H}$.

For the other direction, let $z=X_{H}\left(\mathfrak{z}\right)$ be
an element of $\mathcal{O}_{K}$ for some $\mathfrak{z}\in\left(\mathbb{Z}_{p}\backslash\mathbb{N}_{0}\right)\cap\mathbb{Q}$.
(As for convergence, note that (i) implies that $X_{H}\left(\mathfrak{z}\right)=\lim_{N\rightarrow\infty}X_{H}\left(\left[\mathfrak{z}\right]_{p^{N}}\right)$
converges in $K_{\ell_{n}}$ to $z$ for some $n\geq0$.) Noting the
identity:
\begin{equation}
X_{H}\left(\mathfrak{z}\right)=H_{\left[\mathfrak{z}\right]_{p}}\left(X_{H}\left(\theta_{p}\left(\mathfrak{z}\right)\right)\right)\label{eq:X_H shift identity}
\end{equation}
which follows from rewriting $X_{H}$'s functional equation (\ref{eq:X_H functional equation}),
note that since $\mathfrak{z}$ is in $\left(\mathbb{Z}_{p}\backslash\mathbb{N}_{0}\right)\cap\mathbb{Q}$,
there is an $n\geq0$ so that $\theta_{p}^{\circ n}\left(\mathfrak{z}\right)$'s
sequence of $p$-adic digits is periodic with period $\tau$, being
generated by a string $\mathbf{j}$ of length $\tau$. As such, letting
$\mathfrak{x}=\theta_{p}^{\circ n}\left(\mathfrak{z}\right)$, iterating
(\ref{eq:X_H shift identity}) $\tau$ times gives:
\begin{equation}
X_{H}\left(\mathfrak{x}\right)=H_{\mathbf{j}}\left(X_{H}\left(\theta^{\circ\tau}\left(\mathfrak{x}\right)\right)\right)=H_{\mathbf{j}}\left(X_{H}\left(\mathfrak{x}\right)\right)
\end{equation}
which shows that $X_{H}\left(\mathfrak{x}\right)$ is fixed by $H_{\mathbf{j}}$.
By (ii), using the notion of correctness and the arguments made with
it in \cite{first blog paper}, it follows that for a string $\mathbf{j}$
of length $n$, the extension of $H_{\mathbf{j}}$ to $\mathcal{O}_{K_{v}}$
(which exists, since $v\left(\Lambda\right)\geq1$), $H_{\mathbf{j}}\left(w\right)=w$
occurs if and only if $H_{\mathbf{j}}\left(w\right)=H^{\circ n}\left(w\right)$
and $H^{\circ n}\left(w\right)=w$. This makes $w=X_{H}\left(\mathfrak{x}\right)$
into a periodic point of $H$'s extension to $\mathcal{O}_{K_{v}}$.
Moreover, we have that $X_{H}\left(\mathfrak{z}\right)=z$ lies in
the forward orbit of $X_{H}\left(\mathfrak{x}\right)=w$ under $H$.
By the integrality property of $H$, this forces $w$ to be in $\mathcal{O}_{K}$
since $z$ is, which shows that $z$ and $w$ are both periodic points
of $H$ in $\mathcal{O}_{K}$. $\checkmark$

Q.E.D.

\section{\label{sec:-Adic-Fourier-Analysis}$p$-Adic Fourier Analysis and
the $\ell$-adic Characteristic Function of $X_{H}$}

Here, we introduce the \textbf{Haar measure }on $\mathbb{Z}_{p}$
and \textbf{Pontryagin duality}, the two cornerstones of $p$-adic
Fourier analysis. The Haar measure is the natural generalization of
the Lebesgue measure to locally compact abelian groups such as $\mathbb{Z}_{p}$,
and, like with the Lebesgue integral on $\mathbb{R}$, allows us to
do Fourier analysis with functions on $\mathbb{Z}_{p}$.
\begin{defn}
Letting $p\geq2$, the sets: 
\begin{equation}
k+p^{n}\mathbb{Z}_{p}\overset{\textrm{def}}{=}\left\{ k+p^{n}\mathfrak{z}:\mathfrak{z}\in\mathbb{Z}_{p}\right\} ,\textrm{ }\forall n\geq0,\textrm{ }\forall k\in\left\{ 0,\ldots,p^{n}-1\right\} 
\end{equation}
are both compact and open (``compact-open'') in $\mathbb{Z}_{p}$,
and form a topological base for $\mathbb{Z}_{p}$. As such, they generate
the Borel $\sigma$-algebra on $\mathbb{Z}_{p}$. The \textbf{$p$-adic
Haar measure} $d\mathfrak{z}$ is defined on the Borel $\sigma$-algebra
by the rule: 
\begin{equation}
\int_{k+p^{n}\mathbb{Z}_{p}}d\mathfrak{z}\overset{\textrm{def}}{=}\frac{1}{p^{n}},\textrm{ }\forall n\geq0,\textrm{ }\forall k\in\left\{ 0,\ldots,p^{n}-1\right\} 
\end{equation}
This is the unique translation-invariant probability measure on $\mathbb{Z}_{p}$.
For $r\in\left[1,\infty\right]$, $L^{r}\left(\mathbb{Z}_{p},\mathbb{C}\right)$
is then defined as the Banach space of functions $f:\mathbb{Z}_{p}\rightarrow\mathbb{C}$
under the usual norm: 
\begin{equation}
\left\Vert f\right\Vert _{L^{r}\left(\mathbb{Z}_{p},\mathbb{C}\right)}\overset{\textrm{def}}{=}\begin{cases}
\left(\int_{\mathbb{Z}_{p}}\left|f\left(\mathfrak{z}\right)\right|^{r}d\mathfrak{z}\right)^{1/r} & \textrm{if }r<\infty\\
\textrm{esssup}_{\mathfrak{z}\in\mathbb{Z}_{p}}\left|f\left(\mathfrak{z}\right)\right| & \textrm{if }r=\infty
\end{cases}
\end{equation}
with the standard caveat that each function is actually an equivalence
class of functions whose differences vanish almost everywhere with
respect to $d\mathfrak{z}$.

Next, given $n\geq0$ and $k\in\left\{ 0,\ldots,p^{n}-1\right\} $,
we write:
\begin{equation}
\left[\mathfrak{z}\overset{p^{n}}{\equiv}k\right]\overset{\textrm{def}}{=}\begin{cases}
1 & \textrm{if }\mathfrak{z}\in k+p^{n}\mathbb{Z}_{p}\\
0 & \textrm{else}
\end{cases},\textrm{ }\forall\mathfrak{z}\in\mathbb{Z}_{p}
\end{equation}
This function is said to be \textbf{locally constant mod $p^{n}$},
because its value depends only on the value of $\mathfrak{z}$ mod
$p^{n}$. More generally, a function of a $p$-adic integer variable
is said to be \textbf{locally constant }if it is locally constant
mod $p^{n}$ for some $n\geq0$.

A \textbf{$p$-adic} \textbf{Schwartz-Bruhat (SB) function} is a linear
combination of finitely many indicator functions of this type. Every
SB function is locally constant and Haar integrable, and every locally
constant function is SB and Haar integrable. We write $\mathcal{S}\left(\mathbb{Z}_{p},\mathbb{C}\right)$
to denote the space of SB functions. This is an algebra under point-wise
addition and multiplication, and is dense in $L^{r}\left(\mathbb{Z}_{p},\mathbb{C}\right)$
for all $r$. SB functions are to the classical Schwartz functions
of Fourier analysis on euclidean space to what the Haar measure is
to the Lebesgue.

Finally, although we will not need it for this paper, the Haar measure
naturally extends from $\mathbb{Z}_{p}$ to the field $\mathbb{Q}_{p}$.
Clopen balls of the form $\mathfrak{x}+p^{n}\mathbb{Z}_{p}$ for $\mathfrak{x}\in\mathbb{Q}_{p}$
and $n\in\mathbb{Z}$ form a basis for the topology of $\mathbb{Q}_{p}$.
The $p$-adic Haar measure of $\mathfrak{x}+p^{n}\mathbb{Z}_{p}$
is:
\begin{equation}
\int_{\mathfrak{x}+p^{n}\mathbb{Z}_{p}}d\mathfrak{z}\overset{\textrm{def}}{=}p^{-n},\textrm{ }\forall\mathfrak{x}\in\mathbb{Q}_{p},\textrm{ }\forall n\in\mathbb{Z}
\end{equation}
Note that this makes intuitive sense, as $p^{n}\mathbb{Z}_{p}$ becomes
a $p$-adically \emph{smaller }neighborhood of $0$ as $n\rightarrow+\infty$,
while it becomes a $p$-adically \emph{larger }neighborhood of $0$
as $n\rightarrow-\infty$. 
\end{defn}
\begin{rem}
Henceforth, unless stated otherwise, ``almost everywhere'' , ``measure'',
``measurable'', ``integrable'', etc., are all meant with respect
to the Haar measure.
\end{rem}
\begin{rem}
If $f:\mathbb{Z}_{p}\rightarrow\mathbb{C}$ is continuous (such as
if $f$ is locally constant), it can be shown that:
\begin{equation}
\int_{\mathbb{Z}_{p}}f\left(\mathfrak{z}\right)d\mathfrak{z}=\lim_{N\rightarrow\infty}\frac{1}{p^{N}}\sum_{n=0}^{p^{N}-1}f\left(n\right)\label{eq:Haar integrability}
\end{equation}
Moreover, given any $f:\mathbb{Z}_{p}\rightarrow\mathbb{C}$, the
convergence of the limit (\ref{eq:Haar integrability}) then implies
that $f$ is Haar integrable, and that its integral is given by (\ref{eq:Haar integrability}).
See \cite{second blog paper} for details.
\end{rem}
\begin{rem}
Because $\mathbb{Z}_{p}$ is a finite measure space, for all $r\geq1$
and all measurable $f:\mathbb{Z}_{p}\rightarrow\mathbb{C}$, Hölder's
inequality gives: 
\begin{equation}
\left\Vert f\right\Vert _{L^{r}\left(\mathbb{Z}_{p},\mathbb{C}\right)}\leq\left\Vert f\right\Vert _{L^{qr}\left(\mathbb{Z}_{p},\mathbb{C}\right)},\textrm{ }\forall r,q\geq1
\end{equation}
Thus, $f\in L^{r_{1}}$ implies $f\in L^{r_{2}}$ for all $r_{1},r_{2}\in\left[1,\infty\right]$
with $r_{2}\leq r_{1}$. So, $f\in L^{\infty}\left(\mathbb{Z}_{p},\mathbb{C}\right)$
implies $f\in L^{r}\left(\mathbb{Z}_{p},\mathbb{C}\right)$ for all
$r\in\left[1,\infty\right]$.
\end{rem}
We need the following elementary integration formulae:
\begin{prop}[\textbf{Decomposition formula}]
\label{prop:decomposition}Let $n\geq0$. Then, for all measurable
$f:\mathbb{Z}_{p}\rightarrow\mathbb{C}$:
\begin{equation}
\int_{\mathbb{Z}_{p}}f\left(\mathfrak{z}\right)d\mathfrak{z}=\sum_{k=0}^{p^{n}-1}\int_{k+p^{n}\mathbb{Z}_{p}}f\left(\mathfrak{z}\right)d\mathfrak{z}=\sum_{k=0}^{p^{n}-1}\int_{\mathbb{Z}_{p}}\left[\mathfrak{z}\overset{p^{n}}{\equiv}k\right]f\left(\mathfrak{z}\right)d\mathfrak{z}\label{eq:Change of variables-1}
\end{equation}
as well as the change of variables formula:
\begin{equation}
\int_{a\mathbb{Z}_{p}+b}f\left(\mathfrak{z}\right)d\mathfrak{z}=\left|a\right|_{p}\int_{\mathbb{Z}_{p}}f\left(a\mathfrak{z}+b\right)d\mathfrak{z},\textrm{ }\forall a,b\in\mathbb{Z}_{p}
\end{equation}
Using the change of variables formula, we also have:
\begin{equation}
\int_{\mathbb{Z}_{p}}f\left(\mathfrak{z}\right)d\mathfrak{z}=\frac{1}{p^{n}}\sum_{k=0}^{p^{n}-1}\int_{\mathbb{Z}_{p}}f\left(p^{n}\mathfrak{z}+k\right)d\mathfrak{z}\label{eq:change of variables and decomposition}
\end{equation}
\end{prop}
Next, we have to introduce\textbf{ Pontryagin duality}.\textbf{ }This
is the formalism which allows one to perform Fourier analysis for
functions out of a \textbf{locally compact abelian group} (LCAG).
LCAGs are abelian groups that are equipped with a topological structure
with respect to which the operations of addition and negation are
continuous. $\mathbb{R}$, $\mathbb{C}$, $\mathbb{Z}_{p}$, $\mathbb{Z}$,
$\mathbb{Q}_{p}$, and $\hat{\mathbb{Z}}_{p}$ are all examples of
LCAGs, though for our purposes, we will only need to work with $\mathbb{Z}_{p}$
and $\hat{\mathbb{Z}}_{p}$. We borrow the ``frequency'' terminology
from Tao's blog notes on Fourier Analysis \cite{Tao Fourier Transform Blog Post}.
\begin{defn}
Given an LCAG $G$, written additively, a \textbf{(unitary) character}
of $G$ a continuous group homomorphism $\chi:\left(G,+\right)\rightarrow\left(\mathbb{C}\backslash\left\{ 0\right\} ,\times\right)$
satisfying:
\begin{align}
\chi\left(x+y\right) & =\chi\left(x\right)\chi\left(y\right),\textrm{ }\forall x,y\in G\\
\left|\chi\left(z\right)\right| & =1,\textrm{ }\forall z\in G
\end{align}
The \textbf{trivial character }is the unique character satisfying
$\chi\left(x\right)=1$ for all $x\in G$. A \textbf{frequency} is
a continuous group homomorphism $\xi:\left(G,+\right)\rightarrow\left(\mathbb{R}/\mathbb{Z},+\right)$
so that:
\begin{equation}
\xi\left(x+y\right)=\xi\left(x\right)+\xi\left(y\right),\textrm{ }\forall x,y\in G
\end{equation}
The \textbf{trivial frequency }is the unique frequency satisfying
$\xi\left(x\right)=0$ for all $x\in G$. 

Every unitary character $\chi$ can be uniquely written as:
\begin{equation}
\chi\left(x\right)=e^{-2\pi i\xi\left(x\right)},\textrm{ }\forall x\in G
\end{equation}
for some frequency $\xi$. The \textbf{Pontryagin dual }of $G$, denoted
$\hat{G}$, is the set of all unitary characters, made into a group
by pointwise multiplication; we can also define it as the set of all
frequencies, made into a group by pointwise addition. These are the
multiplicative and additive realizations of the same object. \textbf{Pontryagin
duality }is the observation that the pairing:
\begin{equation}
\left(\xi,x\right)\in\hat{G}\times G\mapsto\xi\left(x\right)\in\mathbb{R}/\mathbb{Z}
\end{equation}
is a $\mathbb{Z}$-bilinear map called the \textbf{duality bracket},
and that this shows that $\hat{G}$ is also an LCAG. The \textbf{Fourier
transform }then sends functions $X:G\rightarrow\mathbb{C}$ to functions
$\hat{X}:\hat{G}\rightarrow\mathbb{C}$ defined by:
\begin{equation}
\hat{X}\left(\xi\right)\overset{\textrm{def}}{=}\int_{G}X\left(x\right)e^{-2\pi i\xi\left(x\right)}dx
\end{equation}
where $dx$ is the Haar measure of $G$, provided the integral exists.
The \textbf{Inverse Fourier Transform }sends functions $\hat{X}:\hat{G}\rightarrow\mathbb{C}$
to functions $X:G\rightarrow\mathbb{C}$:
\begin{equation}
X\left(x\right)=\int_{\hat{G}}\hat{X}\left(\xi\right)e^{2\pi i\xi\left(x\right)}d\xi
\end{equation}
where $d\xi$ is the \textbf{dual measure}, the Haar measure on $\hat{G}$,
provided the integral exists.
\end{defn}
In our case, it can be shown that, as an additive group, the Pontryagin
dual $\hat{\mathbb{Z}}_{p}$ is isomorphic to $\mathbb{Z}\left[1/p\right]/\mathbb{Z}$,
the group of all $p$-power denominator rational numbers in $\left[0,1\right)$
with addition mod $1$ as the group operation. The dual measure on
$\hat{\mathbb{Z}}_{p}$ is just the counting measure. The duality
bracket is given explicitly by the \textbf{$p$-adic fractional part},
$\left\{ \cdot\right\} _{p}:\mathbb{Q}_{p}\rightarrow\hat{\mathbb{Z}}_{p}$.
If $p$ is a prime number, this is defined by: 
\begin{equation}
\left\{ \mathfrak{x}\right\} _{p}=\begin{cases}
0 & \textrm{if }n_{0}\geq0\\
\sum_{n=n_{0}}^{-1}d_{n}p^{n} & \textrm{if }n_{0}\leq-1
\end{cases}
\end{equation}
where:
\begin{equation}
\mathfrak{x}=\sum_{n=n_{0}}^{\infty}d_{n}p^{n}
\end{equation}
is the Hensel series representation of $\mathfrak{x}\in\mathbb{Q}_{p}$.
If $p$ is composite, $\left\{ \cdot\right\} _{p}$ is defined by:
\begin{equation}
\left\{ \mathfrak{x}\right\} _{p}\overset{\textrm{def}}{=}\sum_{\ell\mid p}\left\{ \mathfrak{x}\right\} _{\ell}\label{eq:p-adic fractional part}
\end{equation}
where the sum is taken over all prime divisors $\ell$ of $p$.

Thus, for example, given $t\in\hat{\mathbb{Z}}_{p}$ and $\mathfrak{z}\in\mathbb{Z}_{p}$:
\begin{equation}
e^{2\pi i\left\{ t\mathfrak{z}\right\} _{p}}=\begin{cases}
1 & \textrm{if }t=0\\
e^{2\pi it\left[\mathfrak{z}\right]_{\left|t\right|_{p}}} & \textrm{else}
\end{cases}
\end{equation}
where $\left[\mathfrak{z}\right]_{\left|t\right|_{p}}$ is the projection
of $\mathfrak{z}$ mod the power of $p$ in $t$'s denominator. Thus,
if $p$ is odd:
\begin{equation}
e^{2\pi i\left\{ 2\mathfrak{z}/p^{3}\right\} _{p}}=e^{2\pi i\left(2\left[\mathfrak{z}\right]_{p^{3}}/p^{3}\right)}
\end{equation}
Moreover, we have that:
\begin{equation}
e^{2\pi i\left\{ tn\right\} _{p}}=e^{2\pi itn},\textrm{ }\forall n\in\mathbb{Z},\textrm{ }\forall t\in\hat{\mathbb{Z}}_{p}
\end{equation}
and that:
\begin{equation}
e^{2\pi i\left\{ \frac{t}{n}\right\} _{p}}=e^{2\pi it\left[n\right]_{\left|t\right|_{p}}^{-1}},\forall t\in\hat{\mathbb{Z}}_{p}\backslash\left\{ 0\right\} ,\textrm{ }\forall n\in\mathbb{Z}:\gcd\left(n,p\right)=1
\end{equation}
where $\left[n\right]_{\left|t\right|_{p}}^{-1}$ is the unique integer
in $\left\{ 0,\ldots,\left|t\right|_{p}-1\right\} $ whose product
with $n$ is congruent to $1$ mod $\left|t\right|_{p}$.

Given any function $X\in L^{1}\left(\mathbb{Z}_{p},\mathbb{C}\right)$,
the \textbf{Fourier transform }of $X$, denoted $\hat{X}:\hat{\mathbb{Z}}_{p}\rightarrow\mathbb{C}$
is the function defined by the integral:
\begin{equation}
\hat{X}\left(t\right)\overset{\textrm{def}}{=}\int_{\mathbb{Z}_{p}}X\left(\mathfrak{z}\right)e^{-2\pi i\left\{ t\mathfrak{z}\right\} _{p}}d\mathfrak{z}
\end{equation}
Given a function $\hat{Y}:\hat{\mathbb{Z}}_{p}\rightarrow\mathbb{C}$,
the \textbf{Fourier series generated by }$\hat{Y}$ is:
\begin{equation}
\mathfrak{z}\mapsto\sum_{t\in\hat{\mathbb{Z}}_{p}}\hat{Y}\left(t\right)e^{2\pi i\left\{ t\mathfrak{z}\right\} _{p}}
\end{equation}
where we write $\sum_{t\in\hat{\mathbb{Z}}_{p}}\hat{X}\left(t\right)$
to denote the limit:
\begin{equation}
\sum_{t\in\hat{\mathbb{Z}}_{p}}\hat{X}\left(t\right)\overset{\textrm{def}}{=}\lim_{N\rightarrow\infty}\sum_{\left|t\right|_{p}\leq p^{N}}\hat{X}\left(t\right)\overset{\textrm{def}}{=}\lim_{N\rightarrow\infty}\sum_{k=0}^{p^{N}-1}\hat{X}\left(\frac{k}{p^{N}}\right)\label{eq:t sum}
\end{equation}

\begin{rem}
Most of the qualitative results of classical Fourier analysis hold
in the LCAG case, such as the Riemann-Lebesgue lemma (the Fourier
transform maps functions in $L^{1}$ to functions that decay to $0$
) and Parseval's Identity (the Fourier transform is an isometry from
the Hilbert space $L^{2}\left(G,\mathbb{C}\right)$ to the Hilbert
space $L^{2}\left(\hat{G},\mathbb{C}\right)$), and the \textbf{Convolution
Theorem }(the Fourier transform sends pointwise products to convolutions
and vice-versa.) Of greatest importance to us is the \textbf{Parseval-Plancherel
Identity}:
\begin{equation}
\int_{G}f\left(x\right)g\left(x\right)dx=\int_{\hat{G}}\hat{f}\left(\xi\right)\hat{g}\left(-\xi\right)d\xi\label{eq:PPI}
\end{equation}
valid for all $f$ and $g$ for which either the left-hand or right-hand
sides converge. As with the classical Parseval-Plancherel Identity,
we can use (\ref{eq:PPI}) to define measures and distributions by
specifying how they act in terms of frequencies; ex:
\begin{equation}
\int_{G}f\left(x\right)d\mu\left(x\right)=\int_{\hat{G}}\hat{f}\left(\xi\right)\hat{\mu}\left(-\xi\right)d\xi
\end{equation}
where
\begin{equation}
\hat{\mu}\left(\xi\right)=\int_{G}e^{-2\pi i\xi\left(x\right)}d\mu\left(x\right)
\end{equation}
is the \textbf{Fourier-Stieltjes transform }of the measure $d\mu$
on $G$.
\end{rem}
\begin{rem}
The theory of integration and Fourier analysis we sketched above is
extremely robust, in that it can be transferred almost verbatim to
the case of functions $\mathbb{Z}_{p}\rightarrow R$, where $R$ is
any commutative, unital, metrically complete normed ring so that both
$R$ and $R$'s residue field (in the case where $R$ is a non-archimedean
local ring) have characteristic either $0$ or coprime to $p$. See
\cite{second blog paper} for a comprehensive exposition of this material.
\end{rem}
By \textbf{Lemma \ref{lem:extension and measurability of X_H}}, the
condition $\rho_{H,\ell}<1$ guarantees that $X_{H}:\mathbb{Z}_{p}\rightarrow K_{\ell}$
is measurable. As such, we can view $X_{H}$ as a $K_{\ell}$-valued
random variable whose associated sample space is the Borel $\sigma$-algebra
of $\mathbb{Z}_{p}$. Given any measurable subset $A$ of $K_{\ell}$,
the probability $\mathbb{F}\left(X_{H}\in A\right)$ would then be
the Haar probability measure of the subset of $\mathbb{Z}_{p}$ on
which $X_{H}\in A$. 

In viewing $X_{H}$ as a random variables, it is natural to want to
consider their characteristic function.
\begin{defn}
\label{def:characteristic function}Let $K$ be a number field of
dimension $d$ over $\mathbb{Q}$, and let $\ell$ be a place of $K$.
By algebraic number theory, there exists a unique place $q$ of $\mathbb{Q}$
so that $\ell$ lies over $q$. Letting $\mathbb{Q}_{q}$ denote the
completion of $\mathbb{Q}$ with respect to $q$, given a $p$-Hydra
$H:\mathcal{O}_{K}\rightarrow\mathcal{O}_{K}$, the \textbf{$\ell$-adic
characteristic function of $X_{H}$}, denoted $\hat{\mu}_{H,\ell}$,
is the function $\hat{\mu}_{H,\ell}:\hat{\mathcal{O}}_{K_{\ell}}\rightarrow\mathbb{C}$
defined by:
\begin{equation}
\hat{\mu}_{H,\ell}\left(t\right)\overset{\textrm{def}}{=}\int_{\mathbb{Z}_{p}}e^{-2\pi i\left\{ \textrm{Tr}_{K_{\ell}/\mathbb{Q}_{q}}\left(tX_{H}\left(\mathfrak{z}\right)\right)\right\} _{q}}d\mathfrak{z},\textrm{ }\forall t\in\hat{\mathcal{O}}_{K_{\ell}}\label{eq:phi_H}
\end{equation}
Here, when $\ell$ is non-archimedean, $\hat{\mathcal{O}}_{K_{\ell}}$
is defined to be $K_{\ell}/\mathcal{O}_{K_{\ell}}$, the pontryagin
dual of $\mathcal{O}_{K_{\ell}}$, the ring of integers of the $\ell$-adic
completion of $K$. When $\ell$ is archimedean, $\hat{\mathcal{O}}_{K_{\ell}}$
is defined to be $\mathbb{R}$. $\textrm{Tr}_{K_{\ell}/\mathbb{Q}_{q}}$
is the field trace of the extension $K_{\ell}/\mathbb{Q}_{q}$, and
$\left\{ \cdot\right\} _{q}$ is defined to be the $q$-adic fractional
part $\left\{ \cdot\right\} _{q}:\mathbb{Q}_{q}\rightarrow\hat{\mathbb{Z}}_{q}\subseteq\mathbb{R}/\mathbb{Z}$
when $\ell$ is non-archimedean and is defined to be the identity
map $\left\{ \cdot\right\} _{q}:\mathbb{R}\rightarrow\mathbb{R}$
when $\ell$ is archimedean. $X_{H}$, of course, is the numen of
$H$. We write $e_{\ell}:K_{\ell}\rightarrow\mathbb{T}$ to denote:
\begin{equation}
e_{\ell}\left(\mathfrak{x}\right)=\exp\left(2\pi i\left\{ \textrm{Tr}_{K_{\ell}/\mathbb{Q}_{q}}\left(\mathfrak{x}\right)\right\} \right)_{q},\textrm{ }\forall\mathfrak{x}\in K_{\ell}
\end{equation}
With this notation, we can write $\hat{\mu}_{H,\ell}$ as:
\begin{equation}
\hat{\mu}_{H,\ell}\left(t\right)=\int_{\mathbb{Z}_{p}}e_{\ell}\left(-tX_{H}\left(\mathfrak{z}\right)\right)d\mathfrak{z}
\end{equation}
\end{defn}
\begin{rem}
\label{rem:existence of phi}Note that the measurability of $X_{H}$
guaranteed by $\rho_{H,\ell}$ ensures that the integral defining
$\hat{\mu}_{H,\ell}$ is well-defined and absolutely convergent with
respect to $t$. 
\end{rem}
\begin{prop}
We have:
\begin{equation}
\hat{\mu}_{H,\ell}\left(t\right)=\frac{1}{p}\sum_{j=0}^{p-1}e_{\ell}\left(-c_{j}t\right)\hat{\mu}_{H,\ell}\left(r_{j}t\right)
\end{equation}
\end{prop}
Using these functional equations, one can show that $\hat{\mu}_{H,\ell}$
is the Fourier transform of the unique self-similar probability measure
$\mu_{H,\ell}$ associated to the IFS generated by the $H_{j}$s.
For this, note that by Pontryagin duality, given a function $f:K_{\ell}\rightarrow\mathbb{C}$,
its Fourier transform $\hat{f}:K_{\ell}\rightarrow\mathbb{C}$ is:
\[
\hat{f}\left(w\right)=\int_{K_{\ell}}f\left(z\right)e_{\ell}\left(-wz\right)dz,\textrm{ }\forall w\in K_{\ell}
\]
where $dz$ is the Haar measure on $K_{\ell}$, with $dz$ being normalized
to be the Haar probability measure on $\mathcal{O}_{K_{\ell}}$ if
$\ell$ is non-archimedean. We also need:
\begin{defn}
Let $W\left(K_{\ell},\mathbb{C}\right)$ denote the set of all functions
$\phi:K_{\ell}\rightarrow\mathbb{C}$ so that:
\begin{equation}
\phi\left(z\right)=\int_{K_{\ell}}\hat{\phi}\left(w\right)e_{\ell}\left(wz\right)dw,\textrm{ }\forall z\in K_{\ell}
\end{equation}
for some $\hat{\phi}\in L^{1}\left(K_{\ell},\mathbb{C}\right)$. This
is the \textbf{Wiener algebra }of functions $\phi:K_{\ell}\rightarrow\mathbb{C}$
with absolutely integrable Fourier transforms, and is a Banach algebra
with respect to pointwise multiplication under the norm:
\begin{equation}
\left\Vert \phi\right\Vert _{W\left(K_{\ell},\mathbb{C}\right)}\overset{\textrm{def}}{=}\left\Vert \hat{\phi}\right\Vert _{L^{1}\left(K_{\ell},\mathbb{C}\right)}\overset{\textrm{def}}{=}\int_{K_{\ell}}\left|\hat{\phi}\left(w\right)\right|dw
\end{equation}
When $\ell$ is a non-archimedean place, we write $W\left(\mathcal{O}_{K_{\ell}},\mathbb{C}\right)$
to denote the subspace of $W\left(K_{\ell},\mathbb{C}\right)$ consisting
of all functions in $W\left(K_{\ell},\mathbb{C}\right)$ whose supports
lie in $\mathcal{O}_{K_{\ell}}$. 
\end{defn}
The connection between self-similar measures (\cite{sahlsten survey}
is a good survey article on the subject) and our formalism is then
as follows:
\begin{prop}
\label{prop:self-similarity of d mu_H}Let $d\mu_{H,\ell}$ denote
the continuous linear functional $W\left(K_{\ell},\mathbb{C}\right)\rightarrow\mathbb{C}$
defined by:
\begin{equation}
\int_{K_{\ell}}f\left(z\right)d\mu_{H,\ell}\left(z\right)\overset{\textrm{def}}{=}\int_{K_{\ell}}\hat{f}\left(t\right)\hat{\mu}_{H,\ell}\left(-t\right)dt,\textrm{ }\forall f\in W\left(K_{\ell},\mathbb{C}\right)
\end{equation}
That is:
\begin{equation}
\hat{\mu}_{H,\ell}\left(t\right)=\int_{K_{\ell}}e_{\ell}\left(-tz\right)d\mu_{H}\left(z\right)
\end{equation}
Then, $d\mu_{H,\ell}$ extends to a probability measure on $K_{\ell}$,
and is the unique self-similar Borel probability measure on $K_{\ell}$
satisfying:
\begin{equation}
d\mu_{H,\ell}=\frac{1}{p}\sum_{j=0}^{p-1}H_{j,*}\left\{ d\mu_{H,\ell}\right\} =\frac{1}{p}\sum_{j=0}^{p-1}d\mu_{H,\ell}\circ H_{j}^{-1}
\end{equation}
i.e.:
\begin{equation}
\int_{K_{\ell}}f\left(z\right)d\mu_{H,\ell}\left(z\right)=\frac{1}{p}\sum_{j=0}^{p-1}\int_{K_{\ell}}f\left(H_{j}\left(z\right)\right)d\mu_{H,\ell}\left(z\right)
\end{equation}
In particular, we have that:
\begin{equation}
\int_{K_{\ell}}f\left(z\right)d\mu_{H,\ell}\left(z\right)=\int_{\mathbb{Z}_{p}}f\left(X_{H}\left(\mathfrak{z}\right)\right)d\mathfrak{z},\textrm{ }\forall f\in L^{1}\left(K_{\ell},\mathbb{C}\right)\label{eq:dmu_H as a pushforward of the haar measure}
\end{equation}
\end{prop}
Proof: Self-similarity is straightforward:
\begin{align*}
\int_{K_{\ell}}f\left(z\right)d\mu_{H}\left(z\right) & =\int_{K_{\ell}}\hat{f}\left(t\right)\hat{\mu}_{H,\ell}\left(-t\right)dt\\
\left(\ref{eq:phi_H functional equation}\right); & =\frac{1}{p}\sum_{j=0}^{p-1}\int_{K_{\ell}}\hat{f}\left(t\right)e_{\ell}\left(c_{j}t\right)\hat{\mu}_{H,\ell}\left(-r_{j}t\right)d\mathbf{t}\\
\left(\textrm{Def. of }d\mu_{H,\ell}\right); & =\frac{1}{p}\sum_{j=0}^{p-1}\int_{K_{\ell}}\hat{f}\left(t\right)e_{\ell}\left(c_{j}t\right)\left(\int_{V_{\ell}}e_{\ell}\left(tr_{j}z\right)d\mu_{H,\ell}\left(z\right)\right)dt\\
\left(\textrm{Fubini}\right); & =\frac{1}{p}\sum_{j=0}^{p-1}\int_{K_{\ell}}\left(\int_{K_{\ell}}\hat{f}\left(t\right)e_{\ell}\left(t\left(r_{j}z+c_{j}\right)\right)dt\right)d\mu_{H,\ell}\left(z\right)\\
\left(f\in W\left(K_{\ell},\mathbb{C}\right)\right); & =\frac{1}{p}\sum_{j=0}^{p-1}\int_{K_{\ell}}f\left(r_{j}z+c_{j}\right)d\mu_{H,\ell}\left(z\right)\\
 & =\frac{1}{p}\sum_{j=0}^{p-1}\int_{K_{\ell}}f\left(H_{j}\left(z\right)\right)d\mu_{H,\ell}\left(z\right)
\end{align*}
As for extendibility, the construction of $d\mu_{H,\ell}$ guarantees
that:
\begin{equation}
\int_{K_{\ell}}f\left(z\right)d\mu_{H,\ell}\left(z\right)=\int_{\mathbb{Z}_{p}}f\left(X_{H}\left(\mathfrak{z}\right)\right)d\mathfrak{z},\textrm{ }\forall f\in W\left(K_{\ell},\mathbb{C}\right)\label{eq:dmu_H as a pushforward of the haar measure-1}
\end{equation}
This follows from the fact that:
\begin{align*}
\int_{K_{\ell}}f\left(z\right)d\mu_{H,\ell}\left(z\right) & =\int_{K_{\ell}}\hat{f}\left(w\right)\hat{\mu}_{H,\ell}\left(w\right)dw\\
\left(\textrm{def. of }\hat{\mu}_{H,\ell}\right); & =\int_{K_{\ell}}\hat{f}\left(w\right)\left(\int_{\mathbb{Z}_{p}}e_{\ell}\left(-wX_{H}\left(\mathfrak{z}\right)\right)d\mathfrak{z}\right)dw\\
\left(\textrm{swap integrals}\right); & \overset{!}{=}\int_{\mathbb{Z}_{p}}\left(\int_{K_{\ell}}\hat{f}\left(w\right)e_{\ell}\left(-wX_{H}\left(\mathfrak{z}\right)\right)dw\right)d\mathfrak{z}\\
\left(f\in W\left(K_{\ell},\mathbb{C}\right)\right); & =\int_{\mathbb{Z}_{p}}f\left(X_{H}\left(\mathfrak{z}\right)\right)d\mathfrak{z}
\end{align*}
where the interchange of integrals at the step marked (!) is justified
by the $L^{1}$ convergence of the $dw$ integral. Because of this,
we see that:
\[
\int_{K_{\ell}}f\left(z\right)d\mu_{H,\ell}\left(z\right)\geq0,\textrm{ }\forall f\in W\left(K_{\ell},\mathbb{C}\right):f\left(z\right)\geq0\textrm{ }\forall z\in\mathbb{C}
\]
and hence, that $d\mu_{H,\ell}$ is a \textbf{positive linear functional}.
From this and the containment of $W\left(K_{\ell},\mathbb{C}\right)$
in the space $C_{0}\left(K_{\ell},\mathbb{C}\right)$ of continuous
functions $K_{\ell}\rightarrow\mathbb{C}$ that decay to $0$ as their
input variable tends to $\infty$ in magnitude in $K_{\ell}$, it
follows by the \textbf{Riesz-Markov-Katukani Representation Theorem}
that $d\mu_{H,\ell}$ extends to a Borel probability measure on $K_{\ell}$.
The self-similarity condition proved above guarantees uniqueness.
The aforementioned extension property then lets us extend (\ref{eq:dmu_H as a pushforward of the haar measure-1})
to all $f\in L^{1}$, thus demonstrating (\ref{eq:dmu_H as a pushforward of the haar measure}).

Q.E.D.

\vphantom{}

In the case where $\ell$ is non-archimedean, we can do more. Since
$q$ is the unique prime of $\mathbb{Q}$ lying under $\ell$, it
follows that $K_{\ell}$ is isomorphic to a finite-degree extension
$F$ of $\mathbb{Q}_{q}$. Moreover, $\left|\cdot\right|_{q}$ has
a unique extension to an absolute value on $F$, and this extension
is equal to the absolute value $\left|\cdot\right|_{\ell}$ on $K_{\ell}$.
Thus, we may speak meaningfully about the $q$-adic absolute values
of elements of $K_{\ell}$. To that end, we have the following identities,
which are useful for computations.
\begin{prop}
\label{prop:cardinality of neighborhoods of 0 in O_K_ell-hat}Let
$q$ lie under $\ell$, and let $\pi_{\ell}$ be a uniformizer of
$\mathcal{O}_{K_{\ell}}$, and let $e_{K,\ell}$ be the ramification
index of $K_{\ell}$ over $\mathbb{Q}_{q}$. Then, in $\hat{\mathcal{O}}_{K_{\ell}}=K_{\ell}/\mathcal{O}_{K_{\ell}}$,
for any $n\in\mathbb{N}_{0}$, we have:
\begin{equation}
\left|\left\{ t\in\hat{\mathcal{O}}_{K_{\ell}}:\left|t\right|_{\ell}\leq q^{n/e_{K,\ell}}\right\} \right|=\left(q^{n/e_{K,\ell}}\right)^{\left[K_{\ell}:\mathbb{Q}_{q}\right]}\label{eq:cardinality of neighborhoods of 0 in O_K_ell-hat}
\end{equation}
\end{prop}
\begin{prop}
\label{prop:adjoint identity}Let $m,n$ be non-negative integers
with $n\geq m$, let $A$ be an abelian group, and consider functions
$\hat{f},\hat{g}:\hat{\mathcal{O}}_{K_{\ell}}\longrightarrow A$.
Then:
\begin{equation}
\sum_{\left|t\right|_{\ell}\leq q^{n/e_{K,\ell}}}\hat{f}\left(\pi_{\ell}^{m}t\right)\hat{g}\left(t\right)=\sum_{\left|t\right|_{\ell}\leq q^{\left(n-m\right)/e_{K,\ell}}}\left(\sum_{\left|t\right|_{\ell}\leq q^{m/e_{K,\ell}}}\hat{f}\left(s+\pi_{\ell}^{-m}t\right)\right)\hat{g}\left(t\right)\label{eq:adjoint identity}
\end{equation}
\end{prop}
Finally, using Fourier inversion, we can express the $\ell$-adic
probability distribution of $X_{H}$ in terms of $\hat{\mu}_{H,\ell}$.
\begin{prop}
\label{prop:Fourier inversion of non-archimedean X_H probability}Suppose
that $\ell$ is a non-archimedean place for which $\rho_{H,\ell}<1$
and $\max_{0\leq j<p}\left\Vert r_{j}\right\Vert _{\ell}\leq1$. Then:

\begin{equation}
\mathbb{P}\left(X_{H}\overset{\pi_{\ell}^{n}}{\equiv}\mathfrak{w}\right)=q^{-n\left[K_{\ell}:\mathbb{Q}_{q}\right]/e_{K,\ell}}\sum_{\left|t\right|_{\ell}\leq q^{n/e_{K,\ell}}}\hat{\mu}_{H,\ell}\left(t\right)e_{\ell}\left(t\mathfrak{w}\right),\textrm{ }\forall\mathfrak{w}\in\pi_{\ell}^{-B_{H,\ell}}\mathcal{O}_{K_{\ell}}\label{eq:Fourier inversion of non-archimedean X_H probability}
\end{equation}
for all $n\in\mathbb{N}_{0}$, where $B_{H,\ell}$ is the integer:
\begin{equation}
B_{H,\ell}\overset{\textrm{def}}{=}e_{K,\ell}\max_{0\leq j<p}\log_{q}\left|c_{j}\right|_{\ell}
\end{equation}
where $\log_{q}x=\left(\ln x\right)/\ln q$.
\end{prop}
\begin{rem}
The probability on the left can also be written as:
\begin{equation}
\mathbb{P}\left(X_{H}\overset{q^{n/e_{K,\ell}}}{\equiv}\mathfrak{w}\right)
\end{equation}
 
\end{rem}
Proof: By (\ref{eq:a.e.  ell-adic bound on X_H in NA case}) from
\textbf{Lemma \ref{lem:extension and measurability of X_H}}, the
hypotheses on the $\ell$-adic norms of the $r_{j}$s implies that
$X_{H}$ is $\ell$-adically bounded by: 
\begin{equation}
\max_{0\leq j<p}\left|c_{j}\right|_{\ell}
\end{equation}
for all $\mathfrak{z}$ for which $X_{H}\left(\mathfrak{z}\right)$
is $\ell$-adically convergent. Since the $\ell$-adic absolute value
takes non-zero values of the form $q^{n/e_{K,\ell}}$ for $n\in\mathbb{Z}$,
it follows that there is some constant $B_{H,\ell}\in\mathbb{Z}$
so that:
\[
\max_{0\leq j<p}\left|c_{j}\right|_{\ell}=q^{B_{H,\ell}/e_{K,\ell}}=\left|\pi_{\ell}^{B_{H,\ell}}\right|_{q}
\]
By character orthogonality, using (\ref{prop:cardinality of neighborhoods of 0 in O_K_ell-hat}),
we have:
\[
\sum_{\left|t\right|_{\ell}\leq q^{n/e_{K,\ell}}}e_{\ell}\left(t\left(\mathfrak{y}-\pi_{\ell}^{B_{H,\ell}}\mathfrak{x}\right)\right)=\left(q^{n/e_{K,\ell}}\right)^{\left[K_{\ell}:\mathbb{Q}_{q}\right]}\left[\pi_{\ell}^{B_{H,\ell}}\mathfrak{x}\overset{q^{n/e_{K,\ell}}}{\equiv}\mathfrak{y}\right]
\]
for all $\mathfrak{x}\in K_{\ell}$ with $\left|\mathfrak{x}\right|_{\ell}\leq\pi_{\ell}^{B_{H,\ell}}$
and all $\mathfrak{y}\in\mathcal{O}_{K_{\ell}}$. The congruence in
the Iverson bracket is well-defined because the bound on $\mathfrak{x}$
implies $\pi_{\ell}^{B_{H,\ell}}\mathfrak{x}$ is in $\mathcal{O}_{K_{\ell}}$.
Consequently, our bound on $X_{H}$ implies: 
\[
\sum_{\left|t\right|_{\ell}\leq q^{n/e_{K,\ell}}}e_{\ell}\left(t\left(\mathfrak{y}-\pi_{\ell}^{B_{H,\ell}}X_{H}\left(\mathfrak{z}\right)\right)\right)=\left(q^{n/e_{K,\ell}}\right)^{\left[K_{\ell}:\mathbb{Q}_{q}\right]}\left[\pi_{\ell}^{B_{H,\ell}}X_{H}\left(\mathfrak{z}\right)\overset{q^{n/e_{K,\ell}}}{\equiv}\mathfrak{y}\right]
\]
for all $\mathfrak{y}\in\mathcal{O}_{K_{\ell}}$ and almost every
$\mathfrak{z}\in\mathbb{Z}_{p}$. Integrating with respect to $\mathfrak{z}$
over $\mathbb{Z}_{p}$, we have that:
\[
\int_{\mathbb{Z}_{p}}\left[\pi_{\ell}^{B_{H,\ell}}X_{H}\left(\mathfrak{z}\right)\overset{q^{n/e_{K,\ell}}}{\equiv}\mathfrak{y}\right]d\mathfrak{z}=\mathbb{P}\left(\pi_{\ell}^{B_{H,\ell}}X_{H}\overset{q^{n/e_{K,\ell}}}{\equiv}\mathfrak{y}\right)
\]
 Thus:
\begin{align*}
\left(q^{n/e_{K,\ell}}\right)^{\left[K_{\ell}:\mathbb{Q}_{q}\right]}\mathbb{P}\left(\pi_{\ell}^{B_{H,\ell}}X_{H}\overset{q^{n/e_{K,\ell}}}{\equiv}\mathfrak{y}\right) & =\sum_{\left|t\right|_{\ell}\leq q^{n/e_{K,\ell}}}\int_{\mathbb{Z}_{p}}e_{\ell}\left(t\left(\mathfrak{y}-\pi_{\ell}^{B_{H,\ell}}X_{H}\left(\mathfrak{z}\right)\right)\right)d\mathfrak{z}\\
 & =\sum_{\left|t\right|_{\ell}\leq q^{n/e_{K,\ell}}}e_{\ell}\left(t\mathfrak{y}\right)\hat{\mu}_{H,\ell}\left(\pi_{\ell}^{B_{H,\ell}}t\right)
\end{align*}
Dividing by the constant on the left-hand side gives:
\begin{equation}
\mathbb{P}\left(\pi_{\ell}^{B_{H,\ell}}X_{H}\overset{q^{n/e_{K,\ell}}}{\equiv}\mathfrak{y}\right)=\left(q^{-n/e_{K,\ell}}\right)^{\left[K_{\ell}:\mathbb{Q}_{q}\right]}\sum_{\left|t\right|_{\ell}\leq q^{n/e_{K,\ell}}}\hat{\mu}_{H,\ell}\left(\pi_{\ell}^{B_{H,\ell}}t\right)e_{\ell}\left(t\mathfrak{y}\right)
\end{equation}

Here:
\begin{align*}
\mathbb{P}\left(\pi_{\ell}^{B_{H,\ell}}X_{H}\overset{q^{n/e_{K,\ell}}}{\equiv}\mathfrak{y}\right) & =\mathbb{P}\left(\left|\pi_{\ell}^{B_{H,\ell}}X_{H}-\mathfrak{y}\right|_{q}\leq q^{-n/e_{K,\ell}}\right)\\
 & =\mathbb{P}\left(\left|X_{H}-\frac{\mathfrak{y}}{\pi_{\ell}^{B_{H,\ell}}}\right|_{q}\leq q^{\left(B_{H,\ell}-n\right)/e_{K,\ell}}\right)\\
 & =\mathbb{P}\left(X_{H}\overset{q^{\left(n-B_{H,\ell}\right)/e_{K,\ell}}}{\equiv}\frac{\mathfrak{y}}{q^{B_{H,\ell}/e_{K,\ell}}}\right)
\end{align*}
Since:
\[
\left\{ \pi_{\ell}^{-B_{H,\ell}}\mathfrak{y}:\mathfrak{y}\in\mathcal{O}_{K_{\ell}}\right\} =\pi_{\ell}^{-B_{H,\ell}}\mathcal{O}_{K_{\ell}}
\]
This gives:
\[
\mathbb{P}\left(X_{H}\overset{q^{\left(n-B_{H,\ell}\right)/e_{K,\ell}}}{\equiv}\mathfrak{w}\right)=\left(q^{-n/e_{K,\ell}}\right)^{\left[K_{\ell}:\mathbb{Q}_{q}\right]}\sum_{\left|t\right|_{\ell}\leq q^{n/e_{K,\ell}}}\hat{\mu}_{H,\ell}\left(\pi_{\ell}^{B_{H,\ell}}t\right)e_{\ell}\left(\pi_{\ell}^{B_{H,\ell}}t\mathfrak{w}\right)
\]
for all $\mathfrak{w}\in\pi_{\ell}^{-B_{H,\ell}}\mathcal{O}_{K_{\ell}}$.
Letting $n\geq B_{H,\ell}$, we can apply \textbf{Proposition \ref{prop:adjoint identity}
}to get:\textbf{ }
\begin{align*}
\mathbb{P}\left(X_{H}\overset{q^{\left(n-B_{H,\ell}\right)/e_{K,\ell}}}{\equiv}\mathfrak{w}\right) & =\left(q^{-n/e_{K,\ell}}\right)^{\left[K_{\ell}:\mathbb{Q}_{q}\right]}\sum_{\left|t\right|_{\ell}\leq q^{\left(n-B_{H,\ell}\right)/e_{K,\ell}}}\hat{\mu}_{H,\ell}\left(t\right)e_{\ell}\left(t\mathfrak{w}\right)\sum_{\left|s\right|_{q}\leq q^{B_{H,\ell}/e_{K,\ell}}}1\\
 & =q^{-\left(n-B_{H,\ell}\right)\left[K_{\ell}:\mathbb{Q}_{q}\right]/e_{K,\ell}}\sum_{\left|t\right|_{\ell}\leq q^{\left(n-B_{H,\ell}\right)/e_{K,\ell}}}\hat{\mu}_{H,\ell}\left(t\right)e_{\ell}\left(t\mathfrak{w}\right)
\end{align*}
Replacing $n$ with $n+B_{H,\ell}$ then gives the result.

Q.E.D.


\begin{thebibliography}{10}
\bibitem{Automatic sequences}Allouche, Jean-Paul, and Jeffrey Shallit.
\emph{Automatic sequences: theory, applications, generalizations}.
Cambridge university press, 2003.

\bibitem{Dynamical Systems}Devaney, Robert L. \emph{An introduction
to chaotic dynamical systems}. CRC press, 2018.

\bibitem{Folland - harmonic analysis}Folland, Gerald B. \emph{A course
in abstract harmonic analysis}. Vol. 29. CRC press, 2016.

\bibitem{Automorphic Representations}Goldfeld, D., \& Hundley, J.
(2011). \emph{Automorphic Representations and L-Functions for the
General Linear Group: Volume 1} (Vol. 129). Cambridge University Press.

\bibitem{Hungerford}Hungerford, Thomas W. \emph{Algebra}. Vol. 73.
Springer Science \& Business Media, 2012.

\bibitem{Koblitz's book}Koblitz, Neal. \emph{$p$-adic analysis:
A short course on recent work}. Vol. 46. Cambridge University Press,
1980.

\bibitem{Lang - Algebraic Number Theory}Lang, Serge. \emph{Algebraic
number theory}. Vol. 110. Springer Science \& Business Media, 2013.

\bibitem{Robert's Book}Robert, A. M. (2013). \emph{A course in p-adic
analysis} (Vol. 198). Springer Science \& Business Media. Chicago.

\bibitem{van Rooij - Non-Archmedean Functional Analysis}van Rooij,
A.C.M. \emph{Non-Archimedean functional analysis.} Pure and Applied
Math., vol. 51, Marcel Dekker, New York, 1978.athematicae (Proceedings).
Vol. 72. No. 2. North-Holland, 1969.

\bibitem{Ultrametric Calculus}Schikhof, W. (1985). \emph{Ultrametric
Calculus: An Introduction to p-Adic Analysis} (Cambridge Studies in
Advanced Mathematics). Cambridge: Cambridge University Press. doi:10.1017/CBO9780511623844.

\bibitem{Schikhof's Thesis}Schikhof, W. H. ``Non-Archimedean Harmonic
Analysis'', Ph.D. Thesis, pp. 1\textendash 80 (Catholic Univ. of
Nijmegen, The Netherlands, 1967).

\bibitem{my dissertation}Siegel, Maxwell C. \emph{$\left(p,q\right)$-adic
Analysis and the Collatz Conjecture}, Ph.D. Thesis, (University of
Southern California, 2022). <https://arxiv.org/pdf/2412.02902>. Accessed
15 April 2024.

\bibitem{my paper}Siegel, Maxwell C. ``Infinite Series Whose Topology
of Convergence Varies From Point to Point.'' p-Adic Numbers, Ultrametric
Analysis and Applications 15.2 (2023): 133-167.

\bibitem{first blog paper}\textquotedblleft The Collatz Conjecture
\& Non-Archimedean Spectral Theory: Part I \textemdash{} Arithmetic
Dynamical Systems and Non-Archimedean Value Distribution Theory.\textquotedblright{}
P-Adic Num Ultrametr Anal Appl 16, 143\textendash 199 (2024). https://doi.org/10.1134/S207004662402005

\bibitem{second blog paper}Siegel, M.C. ``The Collatz Conjecture
\& Non-Archimedean Spectral Theory - Part II - $\left(p,q\right)$-Adic
Fourier Analysis and Wiener\textquoteright s Tauberian Theorem''.
P-Adic Num Ultrametr Anal Appl 17, 187\textendash 232 (2025). https://doi.org/10.1134/S2070046625020062

\bibitem{mellin integral blog post}Siegel, Maxwell C. ``The Collatz
Conjecture \& Non-Archimedean Spectral Theory \textendash{} Part II
\textendash{} The Correspondence Principle''. Blog post. <https://siegelmaxwellc.wordpress.com/2022/08/10/the-collatz-conjecture-non-archimedean-spectral-theory-part-ii-the-correspondence-principle/>.
Accessed 14 April 2024.

\bibitem{wiener blog post}Siegel, Maxwell C. ``The Collatz Conjecture
\& Non-Archimedean Spectral Theory \textendash{} Part III \textendash{}
(p,q)-adic Fourier Analysis \& Tauberian Spectral Theory ''. Blog
post. <https://siegelmaxwellc.wordpress.com/2022/08/12/the-collatz-conjecture-non-archimedean-spectral-theory-part-iii-a-pq-adic-fourier-analysis/>.
Accessed 14 April 2024.

\bibitem{Episode 1 video}Siegel, Maxwell C. ``Episode 1 - (p,q)-adic
Analysis and the Collatz Conjecture - A Whole New World''. YouTube,
uploaded by M.C. Siegel, 13 December 2023. <https://youtu.be/xRb8q5DR78E>.

\bibitem{Wiener video}Siegel, Maxwell C. ``Episode 5 (Part 2) -
Where No Wiener Has Gone Before - (p,q)-adic Analysis \& the Collatz
Conjecture''. YouTube, uploaded by M.C. Siegel, 5 February 2024.
<https://youtu.be/nnq-Hu1HJ5Q>.

\bibitem{sahlsten survey}Sahlsten, T. (2025). ``Fourier Transforms
and Iterated Function Systems''. In: Barral, J., Batakis, A., Seuret,
S. (eds) Recent Developments in Fractals and Related Fields. FARF
4 2022. Trends in Mathematics. Birkhäuser, Cham. https://doi.org/10.1007/978-3-031-80453-3\_12

\bibitem{Silverman}Silverman, Joseph H. \emph{The arithmetic of dynamical
system}s. Vol. 241. Springer Science \& Business Media, 2007.

\bibitem{Strogatz}Strogatz, Steven H. \emph{Nonlinear dynamics and
chaos with student solutions manual: With applications to physics,
biology, chemistry, and engineering}. CRC press, 2018.

\bibitem{Tao Fourier Transform Blog Post}Tao, Terence. (2009) ``245C,
Notes 2: The Fourier transform''. <https://terrytao.wordpress.com/2009/04/06/the-fourier-transform/>.
Accessed 12 January 2024.

\bibitem{Vladimirov - the big paper about complex-valued distributions over the p-adics}Vladimirov,
Vasilii S. ``Generalized functions over the field of $p$-adic numbers.''
Russian Mathematical Surveys 43.5 (1988): 19.

\bibitem{Wirsching's book on 3n+1}Wirsching, Günther J. \emph{The
dynamical system generated by the $3n+1$ function}, Lecture Notes
in Mathematics, vol. 1681, Springer-Verlag, Berlin, 1998. MR 1612686.
\end{thebibliography}
\end{document}